\documentclass[onefignum,onetabnum]{siamart220329}

\DeclareSymbolFont{CMlargesymbols}{OMX}{cmex}{m}{n} 
\DeclareMathDelimiter{(}{\mathopen} {operators}{"28}{CMlargesymbols}{"00}
\DeclareMathDelimiter{)}{\mathclose}{operators}{"29}{CMlargesymbols}{"01}
\DeclareMathAlphabet\mathcal{OMS}{cmsy}{m}{n} 
\SetMathAlphabet\mathcal{bold}{OMS}{cmsy}{b}{n} 

\usepackage{listings}
\usepackage{lipsum}
\usepackage{amsfonts}
\usepackage{graphicx}
\usepackage{epstopdf}
\usepackage{algorithmic}
\ifpdf
  \DeclareGraphicsExtensions{.eps,.pdf,.png,.jpg}
\else
  \DeclareGraphicsExtensions{.eps}
\fi

\usepackage[numbers,sort&compress]{natbib}

\newsiamremark{remark}{Remark}
\newsiamremark{hypothesis}{Hypothesis}
\crefname{hypothesis}{Hypothesis}{Hypotheses}
\newsiamthm{claim}{Claim}

\headers{High-Order Fast Direct Solver on Triangles}{G. Zavalani}

\title{A High-Order Fast Direct Solver for Surface PDEs on Triangles
\thanks{Submitted to the editors DATE.
}
}

\author{Gentian Zavalani  \thanks{Technische Universit{\"a}t Dresden, %
  Institute of Numerical Mathematics, %
  01062 Dresden, Germany. \\
Email: \url{gentian.zavalani@tu-dresden.de}
}
}

\usepackage{amsopn}

\ifpdf
\hypersetup{
  pdftitle={Fast spectral methods on triangular and deforming surfaces},
  pdfauthor={G. Zavalani}
}
\fi

\usepackage{multirow} 
\usepackage{mathptmx}
\usepackage[mathcal]{eucal}
\usepackage{nccmath, amssymb, mathtools}
\usepackage[utf8]{inputenc}
\usepackage[T1]{fontenc}
\usepackage{float}
\usepackage{placeins}
\usepackage{overpic}
\usepackage{mathtools}
\usepackage{tgpagella}
\usepackage{multirow}
\usepackage{float}
\usepackage[parfill]{parskip}

\usepackage{graphicx}
\usepackage{lettrine}
\usepackage{enumitem}
\usepackage{mathrsfs}
\usepackage{tikz}
\usetikzlibrary{matrix}
\usepackage{tikz-cd}
\usepackage{booktabs}
\usepackage{caption}
\usepackage{subcaption}
 \definecolor{highlightcolor}{rgb}{0.7, 0.85, 1.0} 

\usepackage{microtype}
\usepackage{bm}
\usepackage{enumitem}
\usepackage{url}
\usepackage{xcolor}
\usepackage{hyperref}
\usepackage{bookmark} 

\definecolor{OurRed}{rgb}{0.64, 0.30, 0.30}

\graphicspath{{images/}}

\makeatletter
\def\input@path{{images/}}
\makeatother

\newcommand{\R}{\mathbb{R}}

\definecolor{MathGreen}{rgb}{0,0.45,0.16}

{\bf}{\it}

\makeatletter
\let\cref@override@label@type\@gobbletwo
\makeatother

\begin{document}

\maketitle

\medskip


\begin{abstract}
In this paper, we develop a triangular formulation of the hierarchical Poincar\'e--Steklov (HPS) method for elliptic partial differential equations on surfaces. It admits high-order discretizations on unstructured triangular meshes of complex geometries, with spectral collocation in the orthogonal Dubiner basis on each element. As in quadrilateral HPS, local solution operators and Dirichlet-to-Neumann maps are built and merged hierarchically, yielding a fast direct solver with $\mathcal O(N \log N)$ complexity for repeated solves on meshes with $N$ elements. Since the precomputed operators are reused for every elliptic solve, the solver also accelerates high-order implicit time-stepping schemes on surfaces. Numerical experiments demonstrate high-order convergence under mesh refinement and spectral accuracy under polynomial refinement for static and time-dependent test problems.

\end{abstract}

\section{Introduction}

Partial differential equations (PDEs) on surfaces arise in a wide range of applications, including fluid dynamics, materials science, and biological pattern formation. Their numerical solution requires discretization techniques that accurately capture both the underlying geometry and the associated surface differential operators. Standard approaches include finite difference, finite element, and spectral methods. Finite difference schemes are simple to implement but offer only algebraic convergence, requiring many grid points and substantial memory. Surface finite element methods (SFEMs)~\cite{dziuk2013finite} provide geometric flexibility but rely on low-order polynomial approximations, typically cubic or quartic~\cite{jonsson2017cut}, which can lead to ill-conditioned systems and make high-order accuracy difficult to achieve, despite the availability of efficient solvers such as multigrid~\cite{briggs2000multigrid} and UMFPACK~\cite{davis2004algorithm}. Spectral methods, by contrast, attain very high accuracy with relatively few degrees of freedom by employing high-degree polynomial approximations, and spectral element methods extend this accuracy to flexible geometries~\cite{fornberg1998practical}, though at the cost of potentially dense linear systems. A particularly effective strategy combines high-order spectral discretizations with fast direct solvers based on the hierarchical Poincaré--Steklov (HPS) method~\cite{gillman2014direct, martinsson2013direct, martinsson2019fast}. In this framework, local solution operators and Dirichlet-to-Neumann maps are constructed on each element and merged hierarchically, yielding efficient solvers with near-linear complexity for repeated solves. The reuse of precomputed operators makes such methods especially attractive for implicit time-stepping schemes for surface PDEs. 

To date, the HPS framework has been most effective on quadrilateral surface meshes~\cite{martinsson2013direct, martinsson2019fast, fortunato2020fast,fortunato2022highorder}, as spectral collocation methods are remarkably efficient on tensor product domains because of the inherent structure of the expansions they employ. The main drawback, however, is the limited ability to handle complex geometries. In such settings, triangular elements are generally preferred due to their geometric flexibility. Therefore, it is tempting to try to marry the efficiency of tensor products with the flexibility of triangular geometries. In~\cite{zavalani2025fast} we explored two ways of doing so: a rhombus-based quadrilateralization of the triangulation, and collocation in the Dubiner basis~\cite{dubiner1991spectral} on Chebyshev--Lobatto-type triangular nodes~\cite{isaac2020recursive}. In~\cite{zavalani2025fast}, the latter is described on a single element, and the merge is taken over unchanged from the quadrilateral case. The present paper develops this approach into a self-contained HPS method on fixed surfaces: we formulate the coupling of neighbouring triangles explicitly and construct the hierarchical merge and the solver from it.

The paper is structured as follows. In Section~\ref{sec:IP}, we describe the Dubiner basis, the triangular nodes, and the resulting interpolation and differentiation operators. Section~\ref{sub:local_discret} presents high-order surface parametrizations and the associated metric tensor. In Sections~\ref{sub:local_op} and~\ref{domain_compos}, we give the local collocation operators and the Schur-complement construction of the HPS merge. Section~\ref{sec:performance} discusses computational complexity. In Section~\ref{time_depdn}, we extend the solver to time-dependent surface PDEs, and Section~\ref{section_numerics} presents numerical examples that demonstrate the speed and high-order convergence of the method on both static and time-dependent problems.

Note that we do not try to compare the quadrilateral and the triangle-based HPS methods in terms of performance or computational cost, as such comparisons depend strongly on the problem setting and implementation details.
Moreover, in the present state of the art, the quadrilateral-based HPS is to be preferred over the triangle-based HPS if a structured
mesh can be used to discretize the computational domain, essentially due to the fact that computing derivatives is more expensive with the triangle-based HPS.

\section{Dubiner polynomials on the reference triangle}\label{sec:IP}

On quadrilaterals, the basis and the nodes come from tensor products of Chebyshev polynomials and points~\cite{fortunato2022highorder}; a triangle has no such tensor-product structure. Following~\cite{zavalani2025fast}, we therefore use the Dubiner polynomials~\cite{dubiner1991spectral}, also called Proriol or Koornwinder polynomials~\cite{proriol1957famille,koornwinder1975two}. They form an orthogonal basis on the reference triangle
$\Delta_2 = \{(\xi,\eta)\in\mathbb{R}^2 : \xi \ge 0,\ \eta \ge 0,\ \xi + \eta \le 1\}$
and are obtained from tensor-product Jacobi polynomials on the square by the collapsed-coordinate mapping.

To keep the paper self-contained, we recall their construction. For total polynomial degree $n$, the Dubiner basis functions are
\[
\Phi_{i,j}(\xi,\eta)
  := C_{ij}\, (1-\eta)^i\,
      J_i^{0,0}\!\left( \frac{2\xi}{1-\eta} - 1 \right)
      J_j^{2i+1,0}(2\eta - 1),
\qquad 0 \le i,\, j \le n,\; i+j \le n,
\]
where the normalization constant is
\[
C_{ij} := \sqrt{2(2i+1)(i+j+1)}.
\]
Here $J_m^{\alpha,\beta}$ denotes the standard Jacobi polynomial, normalized by
$J_m^{\alpha,\beta}(1)=\binom{m+\alpha}{m}$, with weight
\[
w(x) = (1-x)^\alpha (1+x)^\beta,
\]
and orthogonality relation
\[
\begin{gathered}
\int_{-1}^1 (1-x)^\alpha (1+x)^\beta\,
  J_m^{\alpha,\beta}(x)\, J_q^{\alpha,\beta}(x)\, dx
  = h_m^{\alpha,\beta}\,\delta_{mq},\\[3pt]
h_m^{\alpha,\beta}
=\frac{2^{\alpha+\beta+1}\Gamma(m+\alpha+1)\Gamma(m+\beta+1)}
{(2m+\alpha+\beta+1)\,m!\,\Gamma(m+\alpha+\beta+1)}.
\end{gathered}
\]
The formula extends continuously to the vertex $(\xi,\eta)=(0,1)$. With the stated $C_{ij}$, the basis is orthonormal in $L^2(\Delta_2)$ and spans the space $\Pi_{2,n}(\Delta_2)$ of polynomials of total degree at most $n$, which has dimension $N_n = \frac{(n+1)(n+2)}{2}$.

Let $\widehat{X}_n := \left\{ \boldsymbol{\xi}_{i,j} = (\xi_i,\eta_j) \;\middle|\; 0 \le i \le n,\; 0 \le j \le n-i \right\}$
denote a set of spectral nodes on $\Delta_2$, generated by
the recursive, parameter-free construction~\cite{isaac2020recursive}.
In our implementation, the construction starts from the second-kind Chebyshev (Chebyshev--Lobatto) points in one dimension, so the nodes on each edge of the triangle are Chebyshev points. For simplicity, the subscript $(i,j)$ is replaced by a single index $m=1,2,\dots,N_n$ through any fixed bijection $m\equiv m(i,j)$.

For any function $u$ defined on $\Delta_2$, we denote by 
$Q_{\Delta_2,n}u \in \Pi_{2,n}(\Delta_2)$ the $n$th-order polynomial interpolant of $u$ on $\Delta_2$, defined by the interpolation conditions
\[
Q_{\Delta_2,n}u(\boldsymbol{\xi}_m) = u_m := u(\boldsymbol{\xi}_m), \qquad m = 1,\dots,N_n.
\]
The interpolant can be expanded in a nodal or a modal representation. The nodal representation is the Lagrange form
\begin{equation*}
Q_{\Delta_2,n}u(\boldsymbol{\xi}) = \sum_{j=1}^{N_{n}} u_j \, \ell_j(\boldsymbol{\xi}),
\qquad \ell_j(\boldsymbol{\xi}_i)=\delta_{ij}.
\end{equation*}
Since there is no closed-form expression for the Lagrange polynomials through an arbitrary set of points on the triangle~\cite{karniadakis2013spectral}, we work with the modal representation
\begin{equation*}
Q_{\Delta_2,n}u(\boldsymbol{\xi})= \sum_{m=1}^{N_{n}} \bar{u}_m \, \Phi_m(\boldsymbol{\xi}),
\end{equation*}
whose coefficients $\bar u_m$ are not point values. Evaluating it at the nodes gives $\mathbf{u}=\mathcal{V}\bar{\mathbf{u}}$ with the generalized Vandermonde matrix $\mathcal{V}_{j,m} = \Phi_m(\boldsymbol{\xi}_j)$, hence
\begin{equation}\label{eq:modal_coeff}
\bar{\mathbf{u}}=\mathcal{V}^{-1}\mathbf{u},\qquad
\bar{u}_m = \sum_{j=1}^{N_n}(\mathcal{V}^{-1})_{mj}\,u_j,
\end{equation}
and therefore
\begin{equation*}
Q_{\Delta_2,n}u(\boldsymbol{\xi})
= \sum_{m=1}^{N_n}\sum_{j=1}^{N_n}(\mathcal{V}^{-1})_{mj}\,
u_j\,\Phi_m(\boldsymbol{\xi}).
\end{equation*}
Derivatives at the nodes follow by differentiating the modal expansion. With the derivative Vandermonde matrices
\[
(\mathcal{V}_\xi)_{j,m} = \partial_\xi \Phi_m(\boldsymbol{\xi}_j),
\qquad
(\mathcal{V}_\eta)_{j,m} = \partial_\eta \Phi_m(\boldsymbol{\xi}_j),
\]
the nodal derivative vectors are $\mathbf{u}_\xi = \mathcal{V}_\xi\bar{\mathbf u}$ and $\mathbf{u}_\eta = \mathcal{V}_\eta\bar{\mathbf u}$, so that by~\eqref{eq:modal_coeff} the differentiation matrices on the reference triangle are
\begin{equation}\label{eq:diff_tri}
\mathcal{D}_\xi = \mathcal{V}_\xi \mathcal{V}^{-1}, \qquad
\mathcal{D}_\eta = \mathcal{V}_\eta \mathcal{V}^{-1},
\qquad
\mathbf{u}_\xi = \mathcal{D}_\xi \mathbf{u}, \quad
\mathbf{u}_\eta = \mathcal{D}_\eta \mathbf{u}.
\end{equation}
Here $\mathcal{D}_\xi, \mathcal{D}_\eta \in \mathbb{R}^{N_n\times N_n}$, whereas in the quadrilateral case they have dimension $(n+1)^2\times(n+1)^2$, owing to the tensor-product structure.

Before concluding this section, we highlight the relation between the Dubiner basis and the nodal Lagrange basis used, for instance, in~\cite{eskilsson2004triangular}. Comparing the double-sum expansion with the Lagrange form shows that
\[
\ell_j(\boldsymbol{\xi})=\sum_{m=1}^{N_n}(\mathcal V^{-1})_{mj}\,\Phi_m(\boldsymbol{\xi}),
\qquad j=1,\dots,N_n ,
\]
so the Lagrange basis is built from the Dubiner basis through the same Vandermonde matrix. On the same nodes the two represent the same interpolant, and in exact arithmetic they lead to identical collocation operators, with identical conditioning and the same $\mathcal O(N_n^3)=\mathcal O(n^6)$ local cost; they differ only in the roundoff incurred in assembly. We keep nodal unknowns and construct the interpolation and differentiation matrices from the Dubiner basis for the reason given above. The geometry, closest-point projection and degree-$n$ element maps, is that of the high-order Lagrange surface elements of~\cite{Demlow09}.

\section{High-order parametric surface approximation}\label{sub:local_discret}
We consider a general elliptic surface PDE on $\Gamma$,
\begin{equation}\label{eq:problem}
\mathcal{L}_{\Gamma} u(\mathbf{x}) = f(\mathbf{x}), \quad \mathbf{x} \in \Gamma,
\end{equation}
where $f$ is a smooth function on $\Gamma$ and $\mathcal{L}_{\Gamma}$ is a variable-coefficient linear second-order elliptic surface operator. If $\Gamma$ is not closed, \eqref{eq:problem} is supplemented by boundary conditions, e.g.\ $u=h$ on $\partial\Gamma$.
To solve~\eqref{eq:problem} numerically we discretize $\mathcal{L}_{\Gamma}$, $f$, and $h$ in a finite-dimensional basis. The error of a discretization of PDEs on curved surfaces depends on two ingredients: the representation of the solution and the representation of the geometry. A high-order scheme is therefore only possible with a high-order approximation of the surface, which we construct as follows.

Let $\Gamma_h=\bigcup_{k=1}^{N}\widehat{\mathcal E}_k$ be a conforming triangulation of $N$ flat triangles of diameter at most $h$, topologically equivalent to $\Gamma$, with affine maps $\mu_k:\Delta_2\to\widehat{\mathcal E}_k$ from the reference triangle. Let $\pi$ denote the closest-point projection onto $\Gamma$, restricted to $\Gamma_h$. We assume $h$ small enough that $\Gamma_h$ lies in the tubular neighbourhood of $\Gamma$ on which $\pi$ is well defined and bijective, so that $\Gamma=\bigcup_k \pi(\widehat{\mathcal E}_k)$. On this basis we define an $n$th-order polynomial parametrization of the surface.

\begin{figure}[t!]
    \centering
    \begin{tikzpicture}
        \node[inner sep=0pt] at (0,0) {\includegraphics[clip,width=1.0\columnwidth]{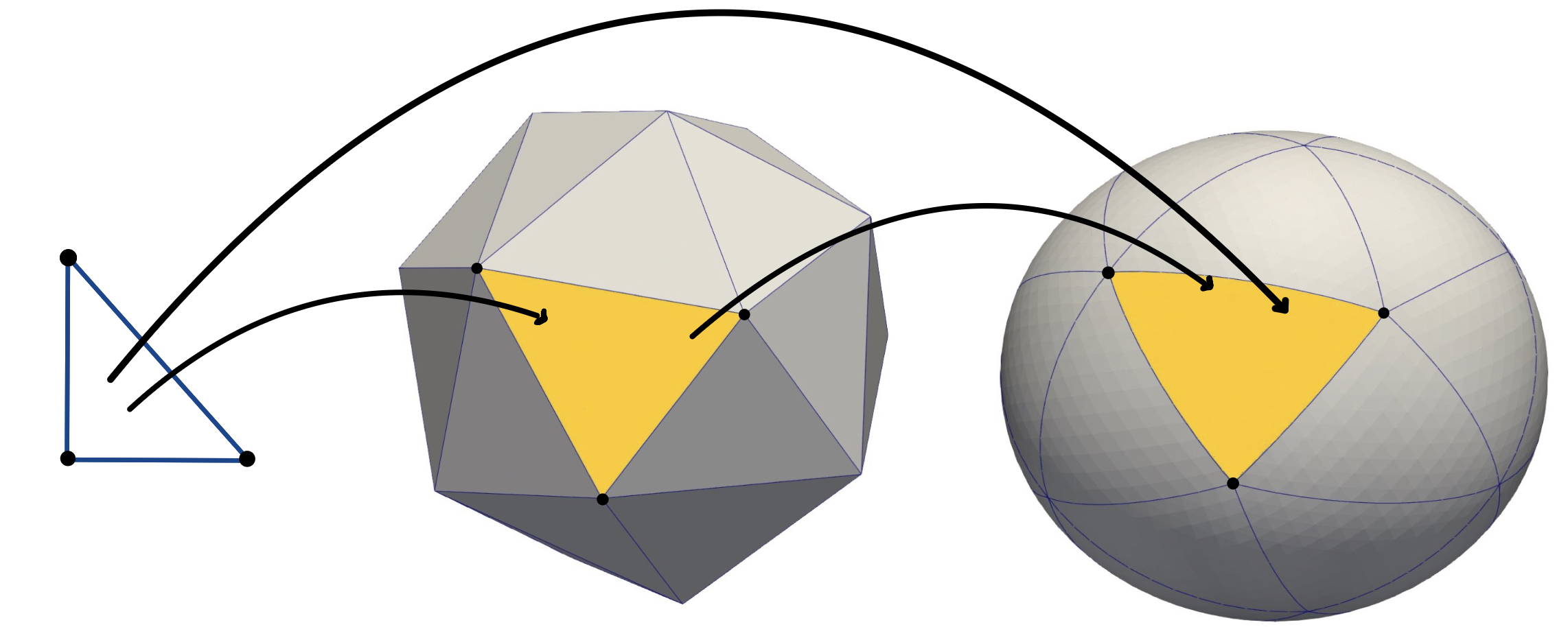}};
        \node[anchor=north west] at (1.0,1.6) {$\pi_k$};
        \node[anchor=north west] at (3.2,-0.1) {$\mathcal{E}_k\,$};
        \node[anchor=north west] at (4.0,2.3) {$\Gamma_{h,n}$};
        \node[anchor=north west] at (-1.5,3.2) {$\phi_{k}:=\pi_{k}\circ\mu_{k}$};
        \node[anchor=north west] at (-2.,-0.1) {$\widehat{\mathcal{E}}_k$};
        \node[anchor=north west] at (-1.8,2.3) {$\Gamma_h$};
        \node[anchor=north west] at (-4.1,0.8) {$\mu_{k}$};
        \node[anchor=north west] at (-6.,-0.6 ) {$\Delta_2$};
    \end{tikzpicture}
  \caption{Construction of a surface parametrization over the reference simplex 
    $\Delta_2$ via closest-point projection from the piecewise affine approximation 
    $\Gamma_{h}$.}
  \label{Figure_quad_app_frame}
\end{figure}

\begin{definition}[Order-$n$ simplex parametrization]\label{def:param}
Let $\phi_k : \Delta_2 \to \mathcal{E}_k\subset\Gamma\subset\mathbb{R}^3$,
where $\mathcal{E}_k:=\phi_k(\Delta_2)$ and
$\phi_k = \pi_k \circ \mu_k$, $k=1,\dots,N$. Here $\pi_k$ is the restriction of $\pi$ to $\widehat{\mathcal{E}}_k$ (see~\Cref{Figure_quad_app_frame}).
We say that the mesh is of order $n$ if each 
element is obtained by polynomial interpolation (in total degree) of $\{\phi_k(\boldsymbol{\xi}_m)\}_{m=1}^{N_n}$, 
sampled at the triangular spectral nodes $\widehat{X}_n$.
\end{definition}

With the Jacobian $D\phi_k(\mathbf{x}):\mathbb{R}^2\to\mathbb{R}^3$, the first fundamental form of $\mathcal E_k$ is the symmetric positive definite matrix
\[
\mathbf{g}(\mathbf{x}) := D\phi_k(\mathbf{x})^T D\phi_k(\mathbf{x}),
\qquad \mathbf{x}\in\Delta_2,
\]
and the volume element is $g(\mathbf{x})=\sqrt{\det\mathbf{g}(\mathbf{x})}$. Appendix~\ref{tang_operator} expresses the tangential operators $\nabla_\Gamma$, $\operatorname{div}_\Gamma$, and $\Delta_\Gamma$ through $\mathbf g$ for a general parametrization.

Since the maps $\phi_k$ can only be evaluated pointwise through the projection, we replace them by their polynomial interpolants. For every triangle we compute $\{\phi_k(\boldsymbol{\xi}_m)\}_{m=1}^{N_n}$ and interpolate the coordinates in total degree $n$,

\begin{equation*}
    Q_{\Delta_2,n}\phi_{k}(\xi,\eta)
    =
    \sum_{m=1}^{N_n}\bar{\boldsymbol\phi}_{k,m}\,\Phi_m(\xi,\eta),
    \qquad
    \bar{\boldsymbol\phi}_{k,m}
    =\sum_{q=1}^{N_n}(\mathcal V^{-1})_{mq}\,\phi_k(\boldsymbol\xi_q).
\end{equation*}
This yields the higher-order approximation $\displaystyle \Gamma_{h,n} := \bigcup_{k=1}^{N}Q_{\Delta_2,n}\phi_k(\Delta_2)$ of $\Gamma$ with surface elements
$\mathcal{K}_{h,n} := \{ Q_{\Delta_2,n}\phi_{k}(\Delta_2) \mid k=1,\dots,N \},$
and the $n$th-order approximation of the volume element
\[
g^n(\mathbf{x}) = \sqrt{\det \left(DQ_{\Delta_2,n}\phi_{k}(\mathbf{x})^T DQ_{\Delta_2,n}\phi_{k}(\mathbf{x})\right)}.
\]

\section{High-order spectral
collocation on a single element}\label{sub:local_op} 
We now describe a spectral collocation method for discretizing~\eqref{eq:problem} on a single surface element $\mathcal{E}_{k}\subset\Gamma$. The discretization is derived directly from the strong form of the problem, and functions on $\mathcal{E}_{k}$ are represented by their values at the same interpolation nodes used for the geometry.

Let $u$ be a function defined on $\mathcal E_k$, with nodal values $u_m:=u(\phi_k(\boldsymbol\xi_m))$. The polynomial interpolant of its pullback has the Dubiner representation
\begin{equation*}\label{eq:u_interp_tri}
Q_{\Delta_2,n}(u\circ\phi_k)(\xi,\eta)
  = \sum_{m=1}^{N_n}\bar u_m\,\Phi_m(\xi,\eta),
  \qquad \bar{\mathbf u}=\mathcal V^{-1}\mathbf u,
  \quad \mathbf u=(u_m)_{m=1}^{N_n}.
\end{equation*}

To discretize $\mathcal{L}_{\Gamma}$ on the element $\mathcal{E}_k$, we first compute the discrete operators on the reference simplex $\Delta_2$ and then map them to $\mathcal{E}_k$ using the numerical coordinate mapping $\mathbf{\phi}_k$. 
Let \(\mathbf{M}[\mathbf{u}] \in \R^{N_n \times N_n}\) denote the diagonal multiplication matrix formed by placing the entries of \(\mathbf{u}\) along the diagonal. Consider the local parametrization $\phi_k(\xi,\eta)
  = \bigl(x_1(\xi,\eta),\, x_2(\xi,\eta),\, x_3(\xi,\eta)\bigr)$.
Using \eqref{eq:diff_tri}, the components of the surface gradient are given by
\begin{equation*}\label{eq:surfgrad_tri}
\mathcal{D}^{\Gamma}_{x_j}
  = \mathbf{M}[\partial_{x_j}\xi]\, \mathcal{D}_\xi
    + \mathbf{M}[\partial_{x_j}\eta]\, \mathcal{D}_\eta,
\qquad j = 1,2,3.
\end{equation*}
For the second-order surface operator\footnote{We identify $\partial_i^{\Gamma}$ with the tangential derivative $\partial_i$ along the surface $\Gamma$.}
\[
\mathcal{L}_\Gamma u
=
\sum_{i,j=1}^3 a_{ij}\,\partial_{x_i}^\Gamma\partial_{x_j}^\Gamma u
+
\sum_{i=1}^3 b_i\,\partial_{x_i}^\Gamma u
+
c\,u,
\]
the discrete operator on the element $\mathcal{E}_k$ is
\begin{equation*}\label{eq:L_tri}
\mathbf{L}_{\mathcal{E}_k}
=
\sum_{i=1}^{3}\sum_{j=1}^{3}
  \mathbf{M}[a_{ij}]\, \mathcal{D}^\Gamma_{x_i}\, \mathcal{D}^\Gamma_{x_j}
\;+\;
\sum_{i=1}^{3} \mathbf{M}[b_i]\mathcal{D}^\Gamma_{x_i}
\;+\;
\mathbf{M}[c].
\end{equation*}

When the principal part is supplied in divergence form,
\[
\mathcal L_\Gamma^{\mathrm{div}}u
=\nabla_\Gamma\!\cdot\!\bigl(A\nabla_\Gamma u\bigr)
+\bm b\cdot\nabla_\Gamma u+cu,
\]
where $A(\mathbf x)$ is uniformly positive definite on the tangent bundle, the same collocation gives the local operator
\[
\mathbf L_{\mathcal E_k}^{\mathrm{div}}
=\sum_{i=1}^{3}\mathcal D_{x_i}^{\Gamma}
 \left(\sum_{j=1}^{3}\mathbf M[A_{ij}]\mathcal D_{x_j}^{\Gamma}\right)
+\sum_{i=1}^{3}\mathbf M[b_i]\mathcal D_{x_i}^{\Gamma}
+\mathbf M[c].
\]

Using this discretization scheme, ~\eqref{eq:problem} can be written as \(N_n\times N_n\) linear system.
\begin{equation}\label{main:disc1}
  \mathbf{L}_{\mathcal{E}_k}\mathbf{u}= \mathbf{f}\;.
\end{equation}
To prepare for the imposition of boundary conditions, we partition the index set \( \{1, \dots, N_n\} \) for a given element \( \mathcal{E}_k \) into interior (\( I_{\operatorname{int}} \)) and boundary (\( I_{\partial} \)) subsets, so that
$\{1, 2, \ldots, N_n \} = I_{\operatorname{int}} \cup I_{\partial},$ with $|I_{\operatorname{int}}| = (n-1)(n-2)/2$ points in the interior and $|I_{\partial}| = 3n$ on the boundary. This partition is used to divide both vectors and matrices into blocks. For instance, given a matrix \( \mathbf{A} \), the submatrix \( \mathbf{A}^{i,b} \) consists of the rows indexed by \( I_{\operatorname{int}} \) and columns indexed by \( I_{\partial} \).~The solution vector is written as
\begin{equation*}
\mathbf{u}^i = \mathbf{u}(I_{\operatorname{int}}) \quad \text{and} \quad \mathbf{u}^b = \mathbf{u}(I_{\partial}).
\end{equation*}

Partitioning the $\mathbf{L}_{\mathcal{E}_{k}}$ and reordering the degrees of freedom in~\eqref{main:disc1} in the order $\{I_{\operatorname{int}}, I_{\partial}\}$ (i.e. interior then boundary) gives a block linear
system,

\begin{equation}\label{main:disc2}
\begin{bmatrix}
\mathbf{L}_{\mathcal{E}_k}^{i,i} & \mathbf{L}_{\mathcal{E}_k}^{i,b} \\
\mathbf{L}_{\mathcal{E}_k}^{b,i} & \mathbf{L}_{\mathcal{E}_k}^{b,b}
\end{bmatrix}
\begin{bmatrix}
\mathbf{u}^i \\
\mathbf{u}^b
\end{bmatrix}
=
\begin{bmatrix}
\mathbf{f}^i \\
\mathbf{f}^b
\end{bmatrix}.
\end{equation}

To impose Dirichlet boundary conditions, we set $\mathbf{u} = \mathbf{h}^b$ on $\partial \mathcal{E}_k$, where $\mathbf{h}^b \in \mathbb{R}^{3n \times 1}$ denotes the vector of boundary values of $h(x)$.  Substituting the prescribed boundary values $\mathbf{u}^b=\mathbf{h}^b$ into the first block row of~\eqref{main:disc2} yields the reduced linear system for the
interior unknowns $\mathbf{u}^i$,
\[
\mathbf{L}_{\mathcal{E}_k}^{i,i} \mathbf{u}^i = \mathbf{f}^i - \mathbf{L}_{\mathcal{E}_k}^{i,b} \mathbf{h}^b.
\]
Define the particular solution and the homogeneous solution operator at the interior nodes by
\[
\bm{v}_{\mathcal{E}_k}^i:=\left(\mathbf{L}_{\mathcal{E}_k}^{i,i}\right)^{-1}\mathbf{f}^i,
\qquad
\mathbf{S}_{\mathcal{E}_k}^i:=-\left(\mathbf{L}_{\mathcal{E}_k}^{i,i}\right)^{-1}\mathbf{L}_{\mathcal{E}_k}^{i,b},
\]
so that
\[
\mathbf{u}^i=\bm{v}_{\mathcal{E}_k}^i+\mathbf{S}_{\mathcal{E}_k}^i\mathbf{h}^b .
\]
We augment these quantities as
\[
\mathbf{S}_{\mathcal{E}_k}:=
\begin{bmatrix}\mathbf{S}_{\mathcal{E}_k}^i\\ I_{3n}\end{bmatrix}
\in\mathbb{R}^{N_n\times 3n},\qquad
\bm{v}_{\mathcal{E}_k}:=
\begin{bmatrix}\bm{v}_{\mathcal{E}_k}^i\\ \mathbf{0}\end{bmatrix},
\]
so that the solution $\mathbf{u}_{\mathcal{E}_k}$ is the sum of the homogeneous solution $\mathbf{w}_{\mathcal{E}_k}:=\mathbf{S}_{\mathcal{E}_k}\mathbf{h}^b$ and the particular solution $\bm{v}_{\mathcal{E}_k}$,
$\mathbf{u}_{\mathcal{E}_k}=\mathbf{w}_{\mathcal{E}_k}+\bm{v}_{\mathcal{E}_k}$.

\section{Domain decomposition methods}\label{domain_compos}

For smooth $u$, the collocation scheme described in the previous section converges rapidly in $n$. The matrix $\mathbf{L}_{\mathcal{E}_k}$ that arises from the discretization has some structure and contains many zeros, but it is still considerably denser than the matrices produced by finite difference~\cite{leveque2007finite} or finite element methods~\cite{lewis2004fundamentals}. One way to reduce this density is to use domain decomposition methods~\cite{mathew2008domain,quarteroni2008numerical}, which also lend themselves well to implementation on parallel computing architectures.

The local solutions $u_k$, each defined only on its element $\mathcal{E}_k$, $k=1,\dots,N$, must fit together into a solution of~\eqref{eq:problem} on the whole surface $\Gamma$. To this end, neighbouring elements $\mathcal{E}_1$ and $\mathcal{E}_2$ are required to satisfy the matching conditions
\begin{equation*}
     u_1 = u_2, \qquad
     \partial_{\mathbf{n}_{b,1}}u_1 = - \partial_{\mathbf{n}_{b,2}} u_2
     \qquad \text{on } \partial \mathcal{E}_{1}\cap\partial \mathcal{E}_{2},
\end{equation*}
continuity of the solution and balance of the outward conormal flux, where $\mathbf{n}_{b,k}$ is the outward unit conormal of $\mathcal{E}_k$: tangent to the surface, normal to the edge.

\begin{remark}\label{rem:conormal}
$\Gamma_{h,n}$ is only $C^0$ across element edges, and the merge works with the conormals of $\Gamma_{h,n}$ itself: each $\mathbf n_{b,k}$ comes from its own element map $Q_{\Delta_2,n}\phi_k$, and the flux balance is imposed with the two as they are. They differ from exactly opposite by the kink angle between the elements, which is $\mathcal O(h^{n})$ for a degree-$n$ surface~\cite{Demlow09}. The exact solution therefore satisfies the discrete interface condition up to $(\mathbf n_{b,1}+\mathbf n_{b,2})\cdot\nabla_\Gamma u=\mathcal O(h^{n}|\nabla_\Gamma u|)$, two orders below the $\mathcal O(h^{n-2})$ observed in Section~\ref{section_numerics}.
\end{remark}

For degree $n\geq3$, each element has $N_n=(n+1)(n+2)/2$ collocation nodes, $3n$ of them on the boundary. In the merge, however, neighbouring elements are coupled on a vertex-free edge grid with $n-1$ first-kind Chebyshev points per edge, so each element has $N_b=3(n-1)$ trace unknowns; the same boundary grid is used in~\cite[Remark~3.1]{fortunato2022highorder}. The matrix $E\in\mathbb R^{3n\times N_b}$ takes edge-grid data to the boundary nodes: it evaluates the degree-$(n-2)$ interpolant on each edge there and averages the two values extrapolated to each vertex. The matrix $R\in\mathbb R^{N_b\times3n}$ goes the other way and interpolates nodal boundary values onto the edge grid. Because the edge grid contains no vertices, no flux equation is imposed at a vertex.

Writing $\mathbf S_{\mathcal E_k}^{\mathrm{nod}}$ for the solution operator from Section~\ref{sub:local_op}, we use
\[
\mathbf S_{\mathcal E_k}:=\mathbf S_{\mathcal E_k}^{\mathrm{nod}}E
\in\mathbb R^{N_n\times N_b},\qquad
\mathbf u_{\mathcal E_k}=\mathbf S_{\mathcal E_k}\mathbf h^b+\bm v_{\mathcal E_k}.
\]
Here $\mathbf h^b$ contains edge-grid data, $\mathcal I_b$ indexes this grid, and $I_{\partial}$ retains its meaning as the nodal boundary index set. Superscripts denote restriction to the stated index sets.

While enforcing continuity of the solution is straightforward, ensuring continuity of the conormal derivative is more delicate. Among the available approaches, an effective strategy is to employ the Poincaré–Steklov operator~\cite{quarteroni2008numerical}, also known as the Dirichlet-to-Neumann (DtN) map. Originally introduced by V.~A.~Steklov, this operator maps prescribed Dirichlet data (i.e., solution values) to the corresponding Neumann data (outward conormal fluxes) on the boundary.
For a given element $\mathcal{E}_k$, the Dirichlet-to-Neumann operator, denoted by $\mathrm{DtN}_{\mathcal{E}_k}$, computes the outward fluxes corresponding to prescribed Dirichlet boundary data. Having constructed the solution operators $\mathbf{S}_{\mathcal{E}_k}$ that solve the PDE locally on each element $\mathcal{E}_k$, for $k = 1, \ldots, N$, the Dirichlet-to-Neumann operator is defined as
\[
\mathrm{DtN}_{\mathcal{E}_k} = \mathcal{D}_{\mathcal{E}_k} \, \mathbf{S}_{\mathcal{E}_k},
\]
where
\[
\mathcal{D}_{\mathcal{E}_k}
=\sum_{j=1}^{3}\mathbf M[\mathbf n_{k,j}]\,
R\,\mathcal D_{x_j}^{\Gamma}(I_{\partial};\cdot)
\in\mathbb R^{N_b\times N_n},
\]
where $\mathbf n_{k,j}$ contains the $j$th component of the outward conormal on the vertex-free edge grid.

For the divergence-form operator, the flux entering the DtN map must use the same principal coefficient. Its boundary operator is
\[
\mathcal F^A_{\mathcal E_k}
=\sum_{i=1}^{3}\sum_{j=1}^{3}
\mathbf M[\mathbf n_{k,i}]\mathbf M[A_{ij}^{\mathrm{tr}}]\,
R\,\mathcal D_{x_j}^{\Gamma}(I_{\partial};\cdot),
\qquad
\mathrm{DtN}^{A}_{\mathcal E_k}
=\mathcal F^A_{\mathcal E_k}\mathbf S_{\mathcal E_k}.
\]
Here $A_{ij}^{\mathrm{tr}}$ denotes coefficient values on the edge grid.
Thus the interface condition balances
$\mathbf n_b\cdot A\nabla_\Gamma u$. For $A=I$, $\mathcal F^A_{\mathcal E_k}=\mathcal D_{\mathcal E_k}$; for $A=aI$, it becomes $\mathbf M[a_b]\mathcal D_{\mathcal E_k}$. The Schur-complement merge is otherwise unchanged.

\subsection{Merging Dirichlet-to-Neumann maps}\label{merge_section}

\begin{figure}[htb]
    \centering
    \begin{overpic}[width=0.4\textwidth]{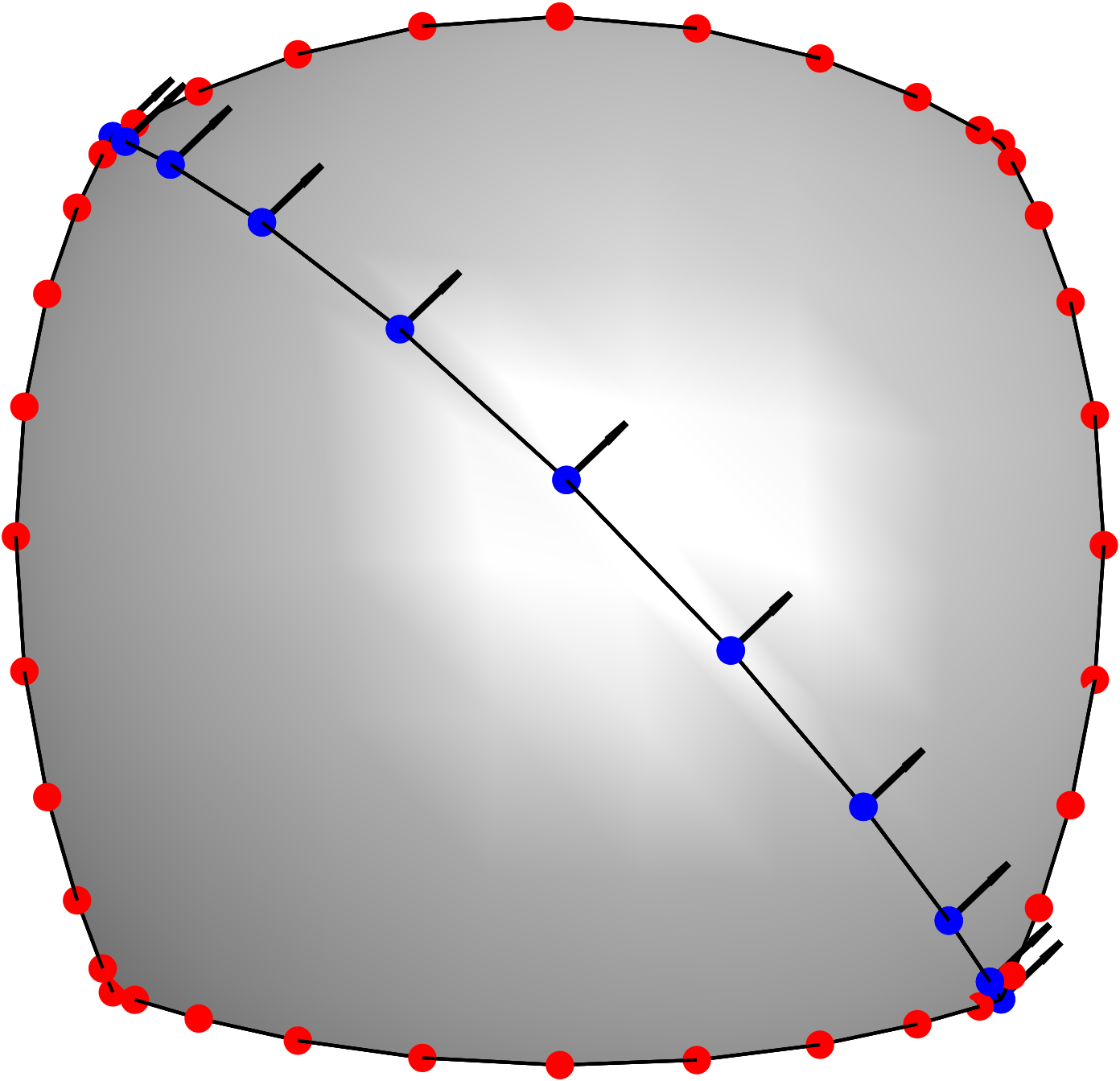}
        \put(32,36){\small$\mathcal{E}_{\bm{1}}$}
        \put(75,45){\small$\mathcal{E}_{\bm{2}}$}
        \put(-20,35){\small$u_1 = h_1$}
        \put(100,59){\small$u_2 = h_2$}
        \put(-23,24){\small$\mathcal{L}_{\Gamma} u_1 = f_1$}
        \put(100,43){\small$\mathcal{L}_{\Gamma}u_2 = f_2$}
        \put(43,53.3){\small$\mathbf{n}_{b}$}
        \put(-24,80){%
          \rotatebox{20}{\small
            $u_1(\mathbf{x}) = u_2(\mathbf{x}) = u_{\text{glue}}(\mathbf{x})$
          }%
        }
        \put(65,-8){%
          \rotatebox{20}{\small
            $\partial_{\mathbf{n}_b} u_1(\mathbf{x}) = -\partial_{\mathbf{n}_b} u_2(\mathbf{x})$
          }%
        }
    \end{overpic}
     \hspace{5em}
    \begin{overpic}[width=0.4\textwidth]{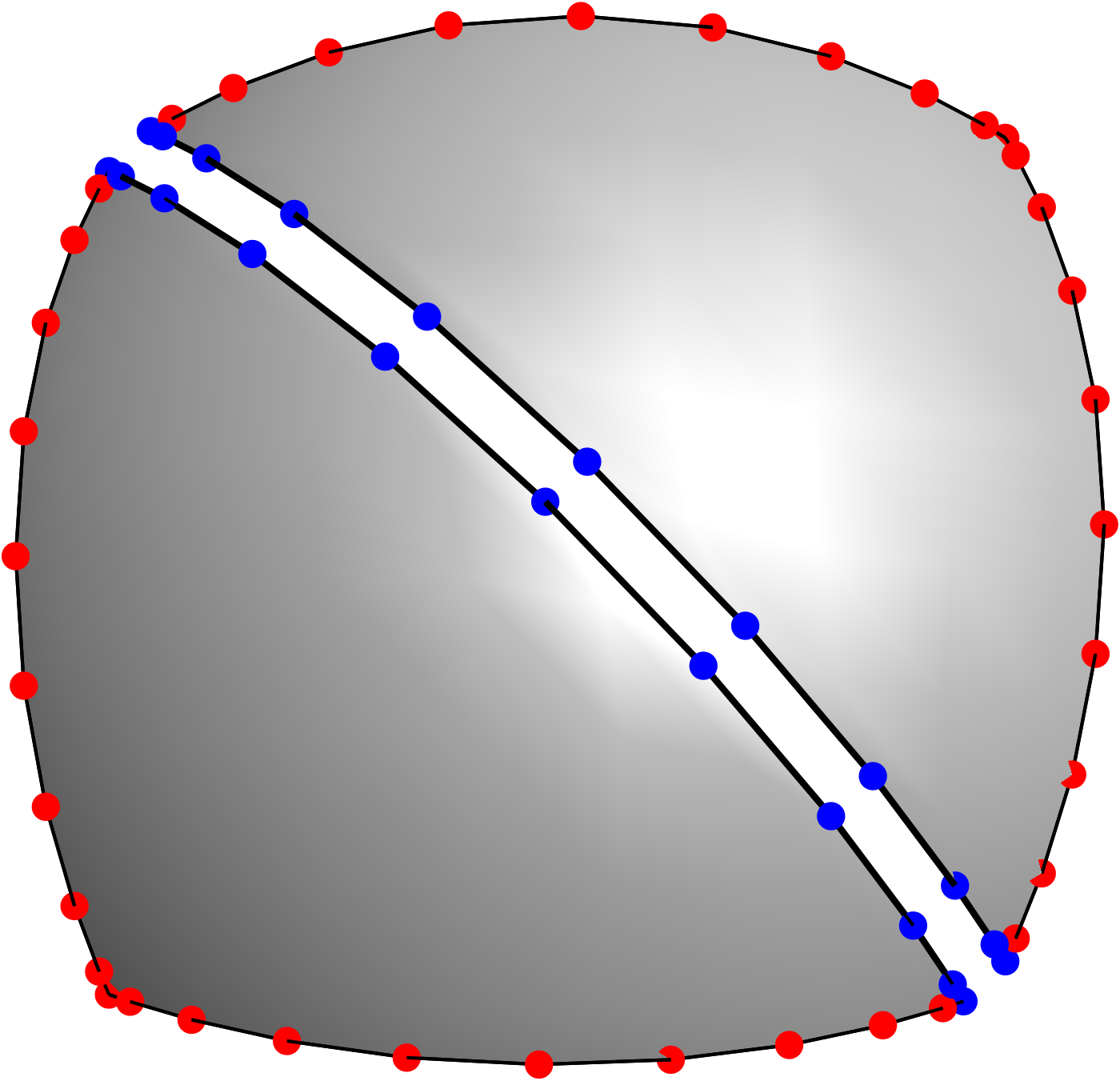}
   \put(-3,16) {\small\color{red}$\mathcal{I}_{b_1}$}
\put(103,66){\small\color{red}$\mathcal{I}_{b_2}$}
\put(40,43){\small\color{blue}$\mathcal{I}_{s_1}$}
\put(55,59){\small\color{blue}$\mathcal{I}_{s_2}$}
    \end{overpic}
\vspace{2em}
\caption{Two surface elements coupled across their shared edge. Left: the solution is continuous and the outward conormal derivatives balance. Right: exterior boundary nodes $\mathcal I_{b_1}$, $\mathcal I_{b_2}$ (red) and aligned interface nodes $\mathcal I_{s_1}$, $\mathcal I_{s_2}$ (blue).}
    \label{fig:tri_two_glued_patches}
\end{figure}

Consider two neighbouring elements $\mathcal{E}_1$ and $\mathcal{E}_2$ with common interface $\mathcal{I}_{12} = \partial \mathcal{E}_1 \cap \partial \mathcal{E}_2$ (\Cref{fig:tri_two_glued_patches}). If the trace $u_{\text{glue}}$ of the solution on $\mathcal I_{12}$ were known, the coupled problem would split into two independent Dirichlet problems,

\begin{equation}\label{eq:glue_1}
\begin{aligned}
\mathcal{L}_{\Gamma} u_k &= f_k && \text{in } \mathcal{E}_k,\\
u_k &= h_k && \text{on } \partial\mathcal{E}_k\setminus \mathcal{I}_{12},\\
u_k &= u_{\text{glue}} && \text{on } \mathcal{I}_{12},
\end{aligned}
\qquad k=1,2,
\end{equation}

each solved by the single-element collocation described in Section~\ref{sub:local_op}. The trace $u_{\text{glue}}$ itself is fixed by the flux condition on $\mathcal I_{12}$. We express it through an interface solution operator $\mathbf S_{\text{glue}}$,
\[
u_{\text{glue}} = \mathbf S_{\text{glue}}
\begin{bmatrix} h_1 \\ h_2 \end{bmatrix},
\]
and now construct $\mathbf S_{\text{glue}}$ from the local Dirichlet-to-Neumann operators of the two elements.

For element $\mathcal E_k$, let $\mathcal I_{s_k}\subset\mathcal I_b$ index the trace grid on the common edge and $\mathcal I_{b_k}:=\mathcal I_b\setminus\mathcal I_{s_k}$ the remaining exterior edges; superscripts $s_k$ and $b_k$ denote row or column restriction to these sets. The two interface node orderings are aligned, and if the edge orientations are opposite the corresponding reversal permutation is applied before coupling. Continuity of the trace then reads $(\mathbf h_1^b)^{s_1}=(\mathbf h_2^b)^{s_2}=: \mathbf u_{\mathrm{glue}}$, so that, with exterior entries first, the boundary data of the two elements are
\[
\mathbf{h}_1^b=\begin{bmatrix}\mathbf{h}_1\\ \mathbf{u}_{\mathrm{glue}}\end{bmatrix},
\qquad
\mathbf{h}_2^b=\begin{bmatrix}\mathbf{h}_2\\ \mathbf{u}_{\mathrm{glue}}\end{bmatrix},
\]
where $\mathbf h_k$ collects the exterior trace values on $\mathcal I_{b_k}$ and $\mathbf u_{\mathrm{glue}}$ is the unknown.

The second interface condition, flux balance, is where the DtN maps enter. With $\mathbf T_k:=\mathrm{DtN}_{\mathcal{E}_k}\in\mathbb R^{N_b\times N_b}$ and $\bm v'_k:=\mathcal D_{\mathcal E_k}\bm v_k$, the outward flux of the particular solution, it reads
\begin{equation}\label{eq:interface_sol}
\bigl(\mathbf T_1\mathbf h_1^b+\bm v'_1\bigr)^{s_1}
+\bigl(\mathbf T_2\mathbf h_2^b+\bm v'_2\bigr)^{s_2}=0,
\end{equation}
a linear system for $\mathbf u_{\mathrm{glue}}$ with matrix $\mathbf A:=\mathbf T_1^{s_1s_1}+\mathbf T_2^{s_2s_2}$. Assuming $\mathbf A$ is invertible, separating the homogeneous and particular contributions gives
\begin{align}
\mathbf w_{\mathrm{glue}}
&=-\mathbf A^{-1}\!\left(\mathbf T_1^{s_1b_1}\mathbf h_1+\mathbf T_2^{s_2b_2}\mathbf h_2\right),\label{eq:interface_homo}\\
\bm v_{\mathrm{glue}}
&=-\mathbf A^{-1}\!\left(\bm v_1'^{s_1}+\bm v_2'^{s_2}\right).\label{eq:interface_partic}
\end{align}
The first of these is the interface solution operator we were after,
\[
\mathbf S_{\mathrm{glue}}
:=-\mathbf A^{-1}\begin{bmatrix}\mathbf T_1^{s_1b_1}&\mathbf T_2^{s_2b_2}\end{bmatrix}.
\]
Substituting the interface trace back into the local DtN maps gives the DtN map of the merged element on its exterior boundary and the corresponding particular flux,
\[
\mathrm{DtN}_{\mathrm{glue}}
=\begin{bmatrix}\mathbf T_1^{b_1b_1}&0\\0&\mathbf T_2^{b_2b_2}\end{bmatrix}
+\begin{bmatrix}\mathbf T_1^{b_1s_1}\\\mathbf T_2^{b_2s_2}\end{bmatrix}\mathbf S_{\mathrm{glue}},
\qquad
\bm v'_{\mathrm{glue}}
=\begin{bmatrix}\bm v_1'^{b_1}\\\bm v_2'^{b_2}\end{bmatrix}
+\begin{bmatrix}\mathbf T_1^{b_1s_1}\\\mathbf T_2^{b_2s_2}\end{bmatrix}\bm v_{\mathrm{glue}},
\]
so the merged element carries exactly the same data as a leaf and can be merged again.

Nothing in these equations depends on element shape: once the local solution operators and DtN maps are defined on matching interface traces, the Schur-complement algebra is the same for triangles as for quadrilaterals, and the triangular geometry and the Dubiner basis enter only through the leaf operators. The merge is therefore applied recursively up the tree. At the root the boundary data, or the closed-surface constraint of Remark~\ref{rem:closed_surface}, are imposed, and a downward pass then recovers the interface traces and the leaf solutions. For a new right-hand side only the particular solutions and their fluxes are recomputed, while the factored interface matrices are kept. The construction is the triangular counterpart of the HPS hierarchy in~\cite{martinsson2019fast,fortunato2022highorder,zavalani2026spectral}.

\begin{remark}\label{rem:closed_surface}
On a connected closed surface, $-\Delta_\Gamma u=f$ requires $\int_\Gamma f\,dS=0$, and we fix $\int_\Gamma u\,dS=0$. As in~\cite[Remark~4.1]{fortunato2022highorder}, a rank-one term regularizes the constant nullspace in the final merge, and we remove the quadrature mean after reconstruction,
\[
\mathbf u\leftarrow\mathbf u-
\frac{\boldsymbol\omega^T\mathbf u}{\boldsymbol\omega^T\mathbf1}\,\mathbf1,
\]
where $\boldsymbol\omega$ holds the elementwise surface quadrature weights.
\end{remark}

\section{Computational asymptotic complexity}\label{sec:performance}

As for quadrilateral HPS~\cite{martinsson2019fast,fortunato2022highorder}, the cost splits into a one-time build-and-factorization stage and a solve stage. The discretization has $N_{\mathrm{dof}}=N N_n$ unknowns before shared traces are identified, where $N$ is the number of elements and $n$ their degree. Building the dense leaf operators costs $\mathcal O(Nn^6)$, since each local collocation matrix has dimension $\mathcal O(n^2)$, and factorizing the merge tree, balanced with separators of the usual square-root size, costs $\mathcal O(N^{3/2}n^3)$, so the one-time cost is $\mathcal O(Nn^6+N^{3/2}n^3)$. A solve with a new right-hand side costs $\mathcal O(n^2N\log N)$ for the hierarchical passes plus $\mathcal O(Nn^4)$ for the dense local particular solves and $\mathcal O(Nn^3)$ for reconstructing the leaf solutions, in total $\mathcal O(n^2N\log N+Nn^4)$. Storing the leaf factors and the hierarchical operators takes $\mathcal O(Nn^4+n^2N\log N)$. At fixed degree these are $\mathcal O(N_{\mathrm{dof}}^{3/2})$ for the factorization and $\mathcal O(N_{\mathrm{dof}}\log N_{\mathrm{dof}})$ per solve.

\section{Time-dependent equations}\label{time_depdn}

Although the triangle-based HPS scheme is designed to solve stationary problems modeled by linear elliptic
PDEs, it can also be useful for accelerating  time-dependent problems modeled by parabolic PDEs. 
We start by considering the reaction-diffusion systems, which are widely regarded as key mechanisms for pattern formation in a variety of contexts, including biological, chemical, physical, and even economic processes. A general reaction--diffusion system describing \(q\) interacting species defined on a closed, smooth surface \(\Gamma \subset \Omega \subset \mathbb{R}^{d+1}\) can be written in the form:
\begin{equation}\label{eq:main_diff}
\frac{\partial \bm{u}}{\partial t} =\nabla_\Gamma \cdot \left( \mathbf{D} \nabla_\Gamma \bm{u} \right)+ \bm{F}(\bm{u}) ,
\end{equation}
where $\bm{u}:=(u_1,\dots,u_s)$, $\mathbf{D}$ is the diffusion tensor, and $\bm{F}$ the reaction kinetics. Spatial discretization with the scheme of the previous sections yields the system of ordinary differential equations
\begin{equation}\label{discret_tdp}
\frac{d\bm{u}}{dt} = L_\Gamma \bm{u} + \bm{F}(\bm{u}).
\end{equation}
The diffusive term $L_\Gamma\bm{u}$ is stiff~\cite{strikwerda2004finite}, so explicit schemes require very small time steps, whereas a fully implicit treatment of $\bm{F}(\bm{u})$ requires its Jacobian, which is typically dense and indefinite~\cite{varga1962matrix}. We therefore use the implicit--explicit backward differentiation formulas (IMEX-BDF)~\cite{ascher1995implicit}, which treat $L_\Gamma\bm{u}$ implicitly and $\bm{F}(\bm{u})$ explicitly.

LLet $\Delta t > 0$ denote the time step and $\bm{u}^n(\mathbf{x}) \approx \bm{u}(\mathbf{x}, n \Delta t)$ the approximate solution at step $n$. The $M$th-order IMEX-BDF scheme reduces each step to the elliptic problem
\begin{equation}
\label{eq:imex_bdf}
A\,\bm{u}^{n+1} = \bm{f}^n,
\qquad
A := I - \omega \Delta t\, L_\Gamma, \qquad
\bm{f}^n := \sum_{i=0}^{M-1} a_i \bm{u}^{n-i} + \Delta t \sum_{i=0}^{M-1} b_i\,\bm{F}(\bm{u}^{n-i}),
\end{equation}
with the coefficients
\[
\begin{aligned}
&\text{IMEX-BDF1:} \; &&\omega = 1, && a = (1), && b = (1), \\
&\text{IMEX-BDF2:} \; &&\omega = \tfrac{2}{3}, && a = \left(\tfrac43, -\tfrac13\right), && b = \left(\tfrac43, -\tfrac23\right), \\
&\text{IMEX-BDF3:} \; &&\omega = \tfrac{6}{11}, && a = \left(\tfrac{18}{11}, -\tfrac{9}{11}, \tfrac{2}{11}\right), && b = \left(\tfrac{18}{11}, -\tfrac{18}{11}, \tfrac{6}{11}\right), \\
&\text{IMEX-BDF4:} \; &&\omega = \tfrac{12}{25}, && a = \left(\tfrac{48}{25}, -\tfrac{36}{25}, \tfrac{16}{25}, -\tfrac{3}{25}\right), && b = \left(\tfrac{48}{25}, -\tfrac{72}{25}, \tfrac{48}{25}, -\tfrac{12}{25}\right).
\end{aligned}
\]
The matrix $A$ is fixed, so its HPS factorization is computed once and each step costs one solve at the price given in Section~\ref{sec:performance}.

\section{Numerical examples}\label{section_numerics}

In this section, we illustrate the accuracy, convergence, and cost of the triangular HPS method on elliptic problems with known solutions and on time-dependent problems on fixed surfaces. The method is implemented in \texttt{pysurfacefun}, an open-source Python package for computing with functions on triangulated surfaces~\cite{pysurfacefun}. All
simulations were preformed in Python~3.14 with NumPy~2.4 on a MacBook Air with an Apple M2 chip and 16\,GB of memory.
\paragraph{Laplace--Beltrami problem on the sphere}
Following the convergence experiments in~\cite{fortunato2022highorder}, we study the Laplace--Beltrami problem on the unit sphere $\Gamma$,
\begin{equation}\label{eq:exp22}
    -\Delta_{\Gamma} u = f \qquad \text{on } \Gamma ,
\end{equation}
with exact solutions taken from the spherical harmonics $Y_\ell^m$, which satisfy
\[
    -\Delta_{\Gamma} Y_\ell^m = \ell(\ell+1)\, Y_\ell^m ,
\]
so that the forcing is $f = \ell(\ell+1)\, Y_\ell^m$. We consider two configurations. On the upper hemisphere we take
\[
    u(\mathbf{x}) = Y_3^2(\mathbf{x}),
\qquad
    h(\mathbf{x}) = 0.25 \sqrt{\frac{105}{\pi}} \, (x_1^2 - x_2^2)\, x_3 \ \text{ on the equator},
\]
and on the closed sphere $u(\mathbf{x}) = Y_{20}^{10}(\mathbf{x})$.

 \begin{figure}[htb]
    \centering
    \begin{subfigure}[t]{0.48\textwidth}
        \centering
        \includegraphics[width=\textwidth]{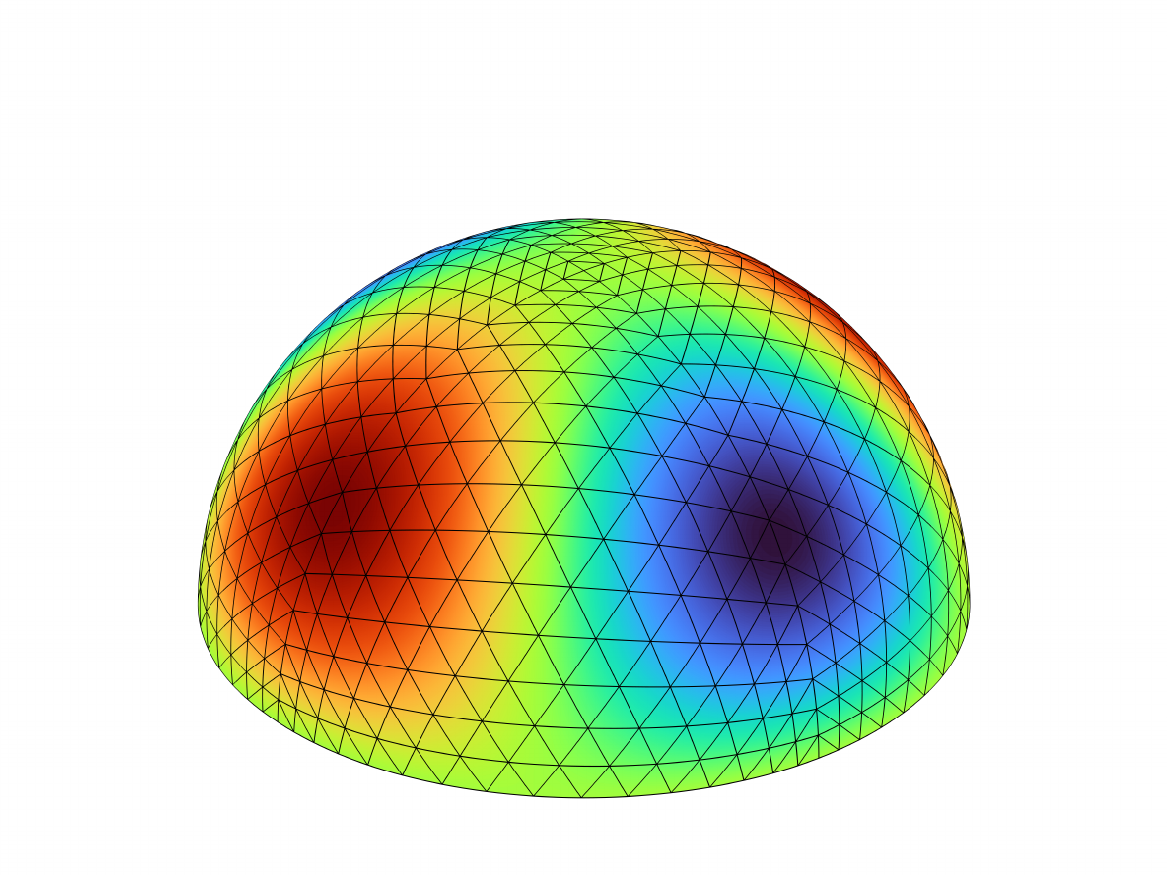}
    \end{subfigure}
    \hfill
    \begin{subfigure}[t]{0.48\textwidth}
        \centering
        \includegraphics[width=\textwidth]{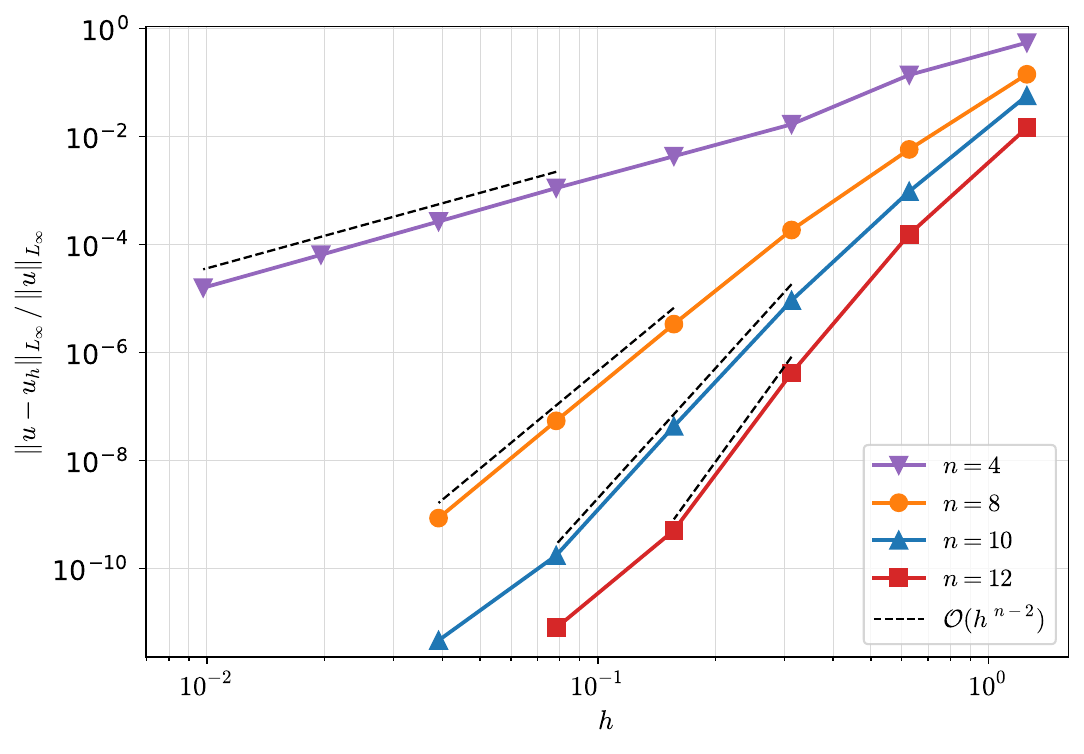}
    \end{subfigure}
    \caption{Laplace--Beltrami equation on the upper hemisphere: computed solution on the triangular mesh (left) and relative $L_\infty$-error under $h$-refinement for $n=4,8,10,12$ with $\mathcal{O}(h^{n-2})$ reference slopes (right).}
    \label{fig:convergence_sphere1}
\end{figure}

\begin{figure}[htb]
    \centering
    \begin{subfigure}[t]{0.5\textwidth}
        \centering
        \includegraphics[width=\textwidth]{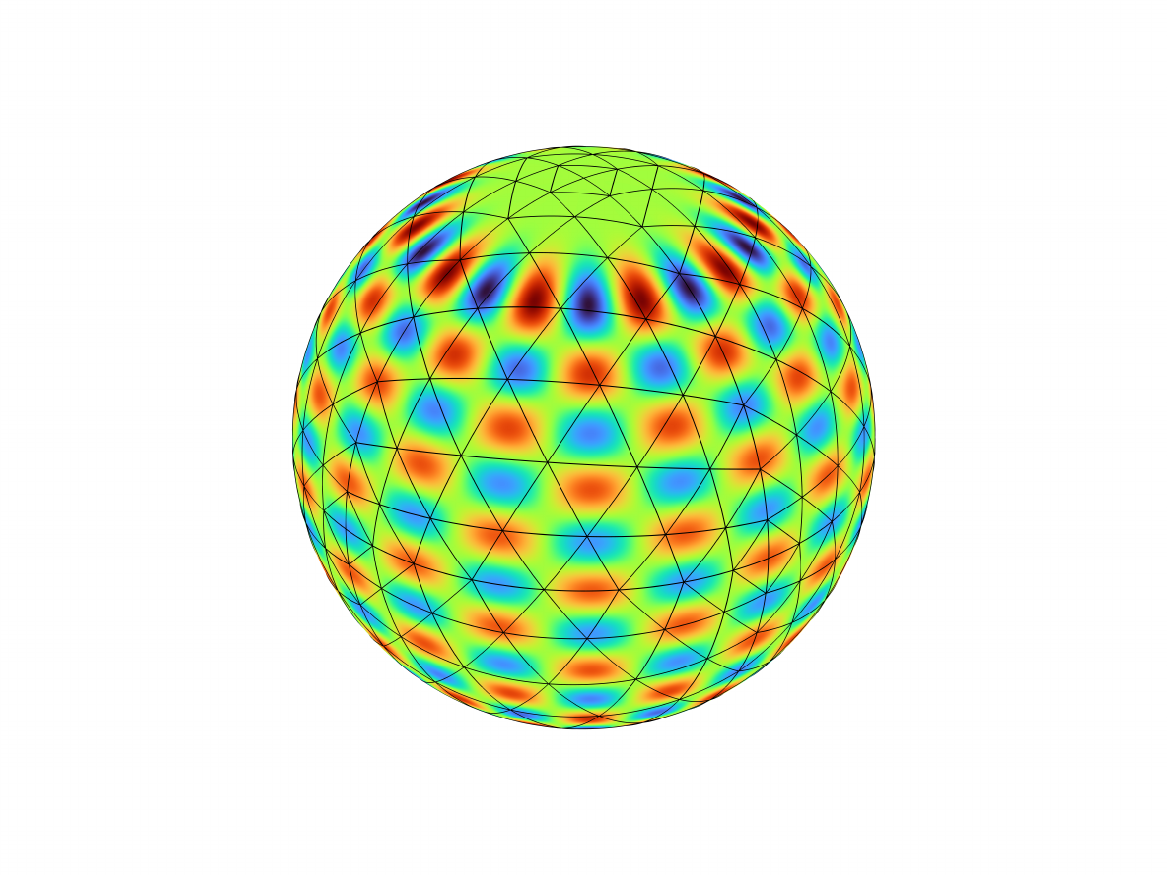}
    \end{subfigure}
    \hfill
    \begin{subfigure}[t]{0.48\textwidth}
        \centering
        \includegraphics[width=\textwidth]{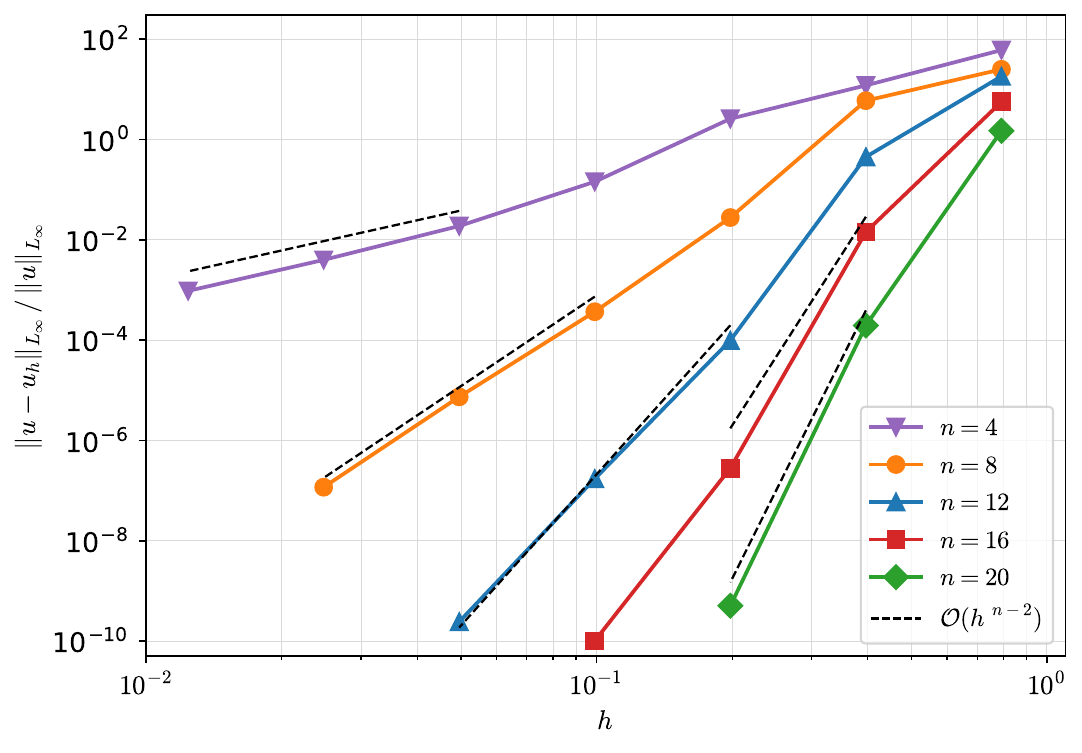}
    \end{subfigure}
\caption{Closed sphere, $u=Y_{20}^{10}$: computed solution (left) and relative $L_\infty$-error under $h$-refinement for $n=4,8,12,16,20$ (right).}
  \label{fig:convergence_sphere}
\end{figure}

\Cref{fig:convergence_sphere1,fig:convergence_sphere} show the computed solutions on the triangulated geometry (left) and the relative $L_\infty$-error under $h$-refinement (right), on the hemisphere and on the closed sphere. As the mesh size $h\to0$, the error decreases at the rate $\mathcal O(h^{n-2})$, indicating high-order convergence;\footnote{The measured slopes over the last refinements above round-off are $2.0$, $6.0$, $7.9$, and $9.7$ on the hemisphere ($n=4,8,10,12$) and $2.1$, $6.0$, and $9.5$ on the closed sphere ($n=4,8,12$).} for $n \geq 16$ the error reaches the rounding floor within two refinements. The rate is two orders below the single-element accuracy because the interface coupling exchanges traces as degree-$(n-2)$ polynomials and differentiates them once in the flux balance. Furthermore, the same degree-$(n-2)$ boundary grid is used in the quadrilateral method of~\cite{fortunato2022highorder}, where the rate is reported as $\mathcal O(h^{n-1})$ over four to five refinements. Refining its sphere test two levels further with our implementation of the quadrilateral method, which reproduces the published errors on the finest meshes there, gives slopes of $2.1$ for $n=4$ and $6.1$ for $n=8$, so we expect both methods to converge at the rate $\mathcal O(h^{n-2})$. An a priori convergence theorem for this class of methods is not yet available~\cite{Demlow09}.

\begin{figure}[htbp]
  \centering
  \includegraphics[width=\textwidth]{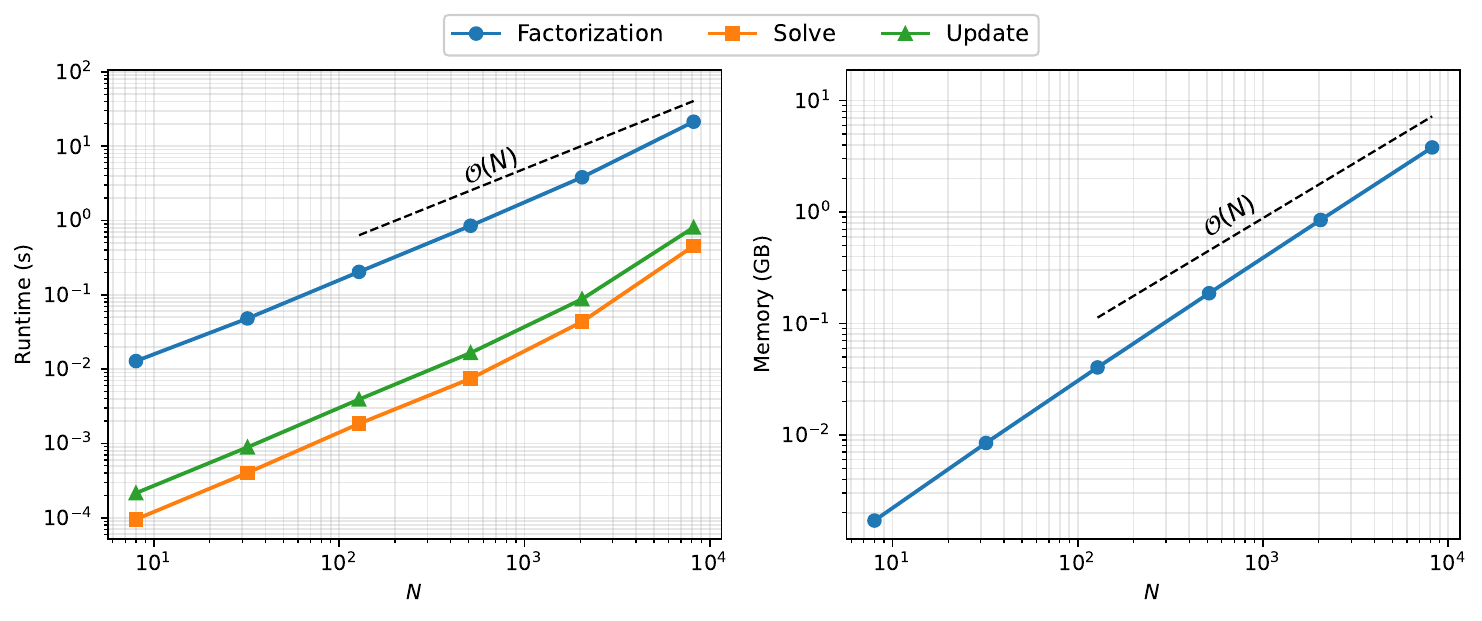}
  \caption{Performance of the triangular HPS method for the Laplace--Beltrami problem on the sphere at degree $n=15$. Left: runtime of the factorization, solve, and right-hand-side update stages as the number of elements $N$ increases. Right: memory required to store the HPS operators.}
  \label{fig:timing_scaling}
\end{figure}

\paragraph{Performance}
To illustrate the cost of the triangular HPS method, we measure the runtime and memory of the factorization and solve stages as the number $N$ of elements grows. We use the closed-sphere Laplace--Beltrami problem with $u=Y_{20}^{10}$ at degree $n=15$, on meshes obtained by projecting uniformly refined octahedra onto the sphere, up to about $10^{6}$ degrees of freedom. \Cref{fig:timing_scaling} (left) shows the runtime of the factorization, of a solve, and of an update of the particular solutions for a new right-hand side followed by a solve. The factorization grows linearly in $N$ over most of the range, since at degree $15$ the leaf term $\mathcal O(Nn^6)$ of Section~\ref{sec:performance} dominates; at the finest mesh it begins to bend upward as the merge term $\mathcal O(N^{3/2}n^3)$ takes over. A solve and an update remain one to two orders of magnitude cheaper than the factorization at every size, so in time stepping, where the factorization is computed once, each step costs only an update and a solve. \Cref{fig:timing_scaling} (right) shows the memory needed to store the HPS operators, which grows linearly in $N$ over the whole range.

\begin{figure}[htbp]
  \centering
  \begin{subfigure}[t]{0.48\textwidth}
    \centering
    \includegraphics[width=\textwidth]{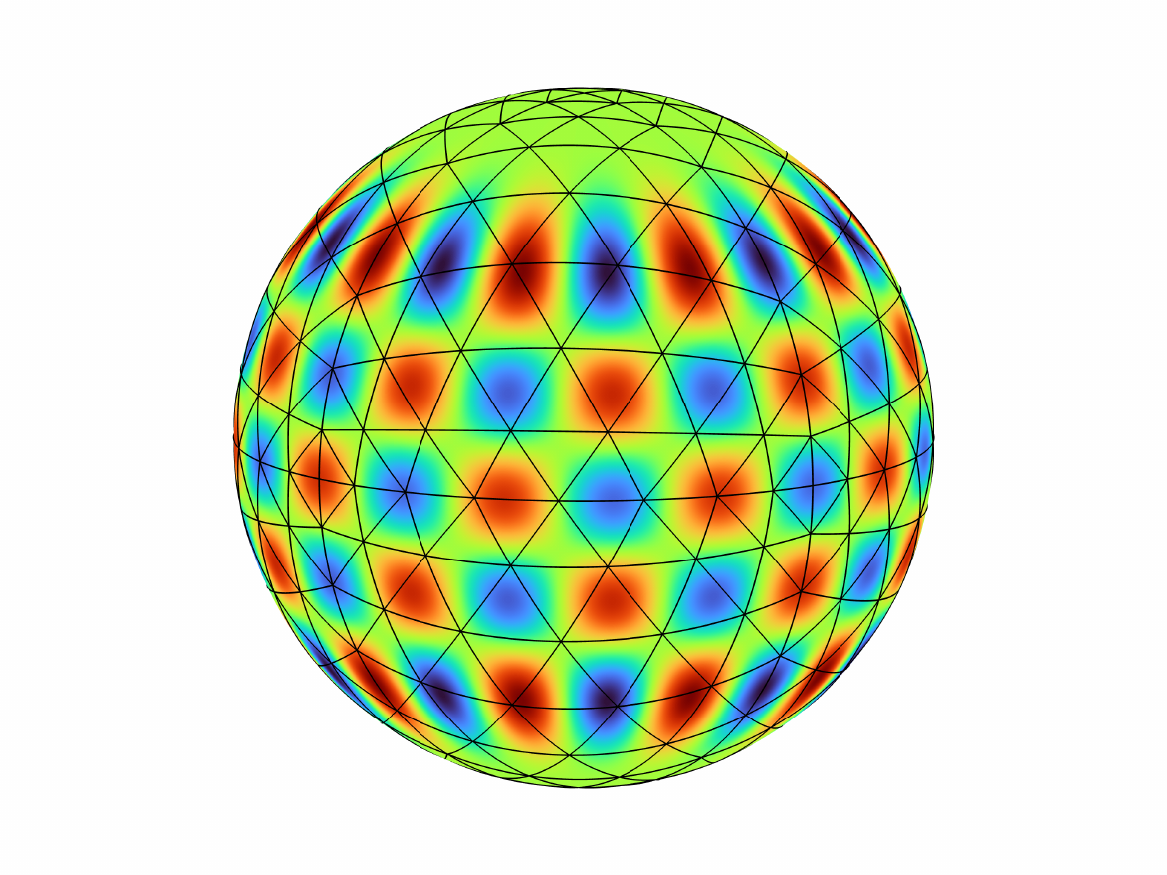}
  \end{subfigure}
  \hfill
  \begin{subfigure}[t]{0.48\textwidth}
    \centering
    \includegraphics[width=\textwidth]{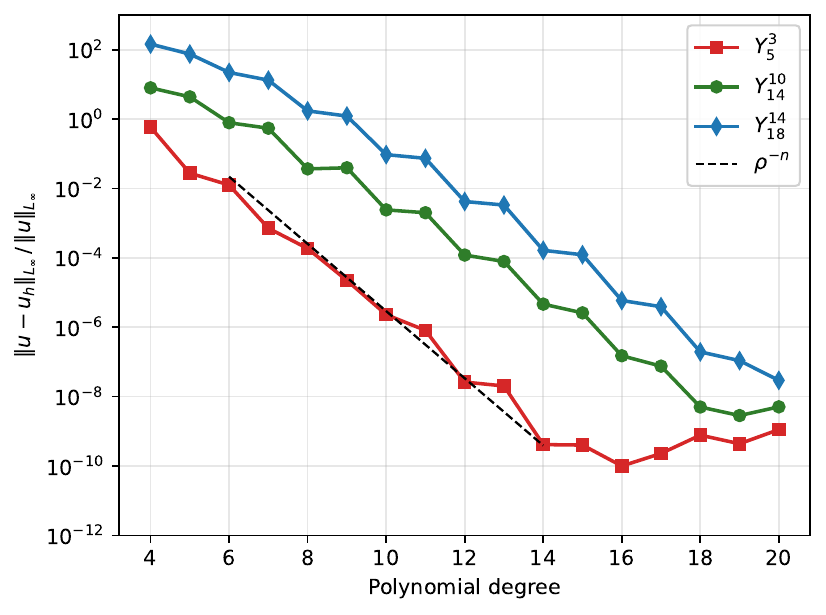}
  \end{subfigure}
  \caption{Surface biharmonic problem on the unit sphere. Left: representative spherical-harmonic solution on the high-order triangular mesh. Right: relative $L_\infty$-error under polynomial refinement for three harmonic eigenfunctions; the dashed line is an exponential reference rate.}
  \label{fig:biharmonic_prefinement_sphere}
\end{figure}

\paragraph{Surface biharmonic operator}
We next consider the surface biharmonic equation on the unit sphere,
\[
(-\Delta_\Gamma)^2u=f ,
\]
with manufactured solutions $u=Y_\ell^m$, for which the eigenrelation above gives $f=\ell^2(\ell+1)^2\,Y_\ell^m$. Introducing $w=-\Delta_\Gamma u$ reduces the fourth-order equation to two sequential Poisson problems,
\begin{equation}\label{eq:biharmonic_split}
-\Delta_\Gamma w=f,
\qquad
-\Delta_\Gamma u=w,
\end{equation}
both on a closed surface and therefore handled as in Remark~\ref{rem:closed_surface}: the surface-quadrature mean is removed from $f$, from the intermediate solution $w$, and from $u$. The split reuses the Laplace--Beltrami factorization twice and requires neither a fourth-order leaf operator nor a fourth-order DtN construction. Figure~\ref{fig:biharmonic_prefinement_sphere} reports a $p$-refinement study on a fixed mesh for three spherical harmonics: the relative $L_\infty$-error decays exponentially with the polynomial degree until it reaches the accuracy floor, so the split retains the spectral convergence of the underlying second-order solver.

\begin{figure}[htbp]
  \centering
  \includegraphics[width=0.58\textwidth]{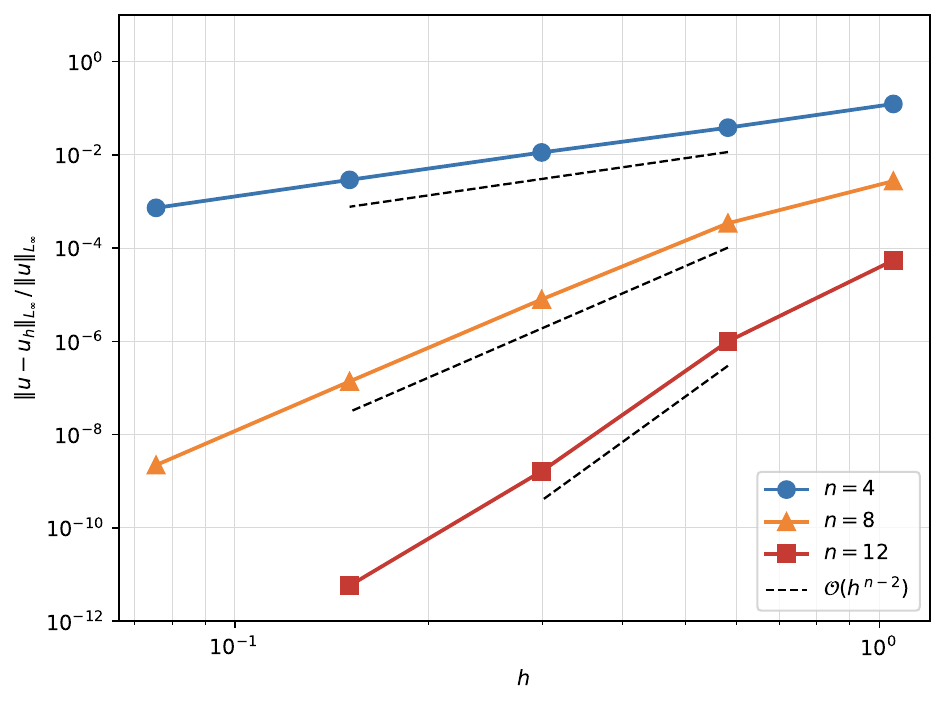}
  \caption{Variable-coefficient divergence-form problem on the unit sphere with $a=2+x_1$, $b=-1$, and exact solution $u=x_1x_2x_3$. Relative $L_\infty$-error under mesh refinement for polynomial degrees $n=4,8,12$, shown with $\mathcal O(h^{n-2})$ reference slopes.}
  \label{fig:divform_sphere}
\end{figure}

\paragraph{Variable-coefficient divergence-form problem on the sphere}
To test the expanded-form discretization of Section~\ref{sub:local_op} with a smooth scalar coefficient, we solve
\[
\nabla_\Gamma\!\cdot(a\nabla_\Gamma u)+bu=f
\qquad\text{on the unit sphere }\Gamma,
\]
with the strictly positive diffusivity $a(\mathbf x)=2+x_1$ and the shift $b=-1$, which makes the operator nonsingular. For the manufactured solution $u=x_1x_2x_3$, which satisfies $\Delta_\Gamma u=-12u$, the forcing is
\[
\begin{aligned}
f
&=a\Delta_\Gamma u+\nabla_\Gamma a\cdot\nabla_\Gamma u+bu\\
&=(2+x_1)(-12x_1x_2x_3)
+x_2x_3(1-3x_1^2)+b\,x_1x_2x_3 ,
\end{aligned}
\]
where $\nabla_\Gamma a\cdot\nabla_\Gamma u$ is evaluated with the surface projector $P_\Gamma=I-\mathbf n_\Gamma\mathbf n_\Gamma^{T}$, $\mathbf n_\Gamma=\mathbf x$. The merge balances the unweighted conormal derivatives; since $a$ is continuous and strictly positive across each interface, this is equivalent to balancing the flux $a\,\partial_{\mathbf n_b}u$ at the matching points. Figure~\ref{fig:divform_sphere} shows the relative $L_\infty$-error under $h$-refinement for $n=4,8,12$. The order increases with the degree and follows the same $\mathcal O(h^{n-2})$ slopes as the constant-coefficient problem, so the smooth variable coefficient costs no accuracy.

\paragraph{Surface diffusion equation on the sphere}
Next, we assess the accuracy and convergence of the triangle-based HPS scheme for a time-dependent problem. As a benchmark with a known analytical solution, we consider isotropic diffusion on the surface of the unit sphere,

\begin{equation*}
    \frac{\partial u}{\partial t} = \Delta_{\Gamma} u \quad \text{on } \Gamma = \{\mathbf{x} \in \mathbb{R}^3 \mid \|\mathbf{x}\|_2 = 1\}.
\end{equation*}
The exact solution in spherical coordinates \( (\vartheta, \varphi) \), where \( \vartheta \in [0,\pi] \) is the polar angle and \( \varphi \in [0, 2\pi) \) is the azimuthal angle, is obtained via expansion in spherical harmonics \( Y_\ell^m \):
\begin{equation*}
    u(t, \vartheta, \varphi) = \sum_{\ell = 0}^{\infty} \sum_{m = -\ell}^{\ell} c_{\ell m}(0) Y_\ell^m(\vartheta, \varphi) e^{- \ell(\ell + 1)t},
\end{equation*}
with modal coefficients determined by the projection of the initial data:
\begin{equation*}
    c_{\ell m}(0) = \int_{\Gamma} (-1)^m Y_\ell^{-m}(\vartheta, \varphi)\, u(0, \vartheta, \varphi)\, dS.
\end{equation*}

To avoid truncating this expansion, we take $u(0,\vartheta,\varphi)=Y_1^0(\vartheta,\varphi)=\sqrt{3/(4\pi)}\cos\vartheta$, with polar angle $\vartheta$, so that the exact solution is $u(t,\vartheta,\varphi)=Y_1^0(\vartheta,\varphi)\,e^{-2t}$.
We integrate to $t=1$ with IMEX-BDF4 and $\Delta t=10^{-3}$, starting the multistep scheme from the exact solution at $t_j=j\Delta t$, $j=0,\dots,3$.

\begin{figure}[htbp]
  \centering
  \begin{subfigure}[b]{0.45\textwidth}
    \includegraphics[width=\textwidth]{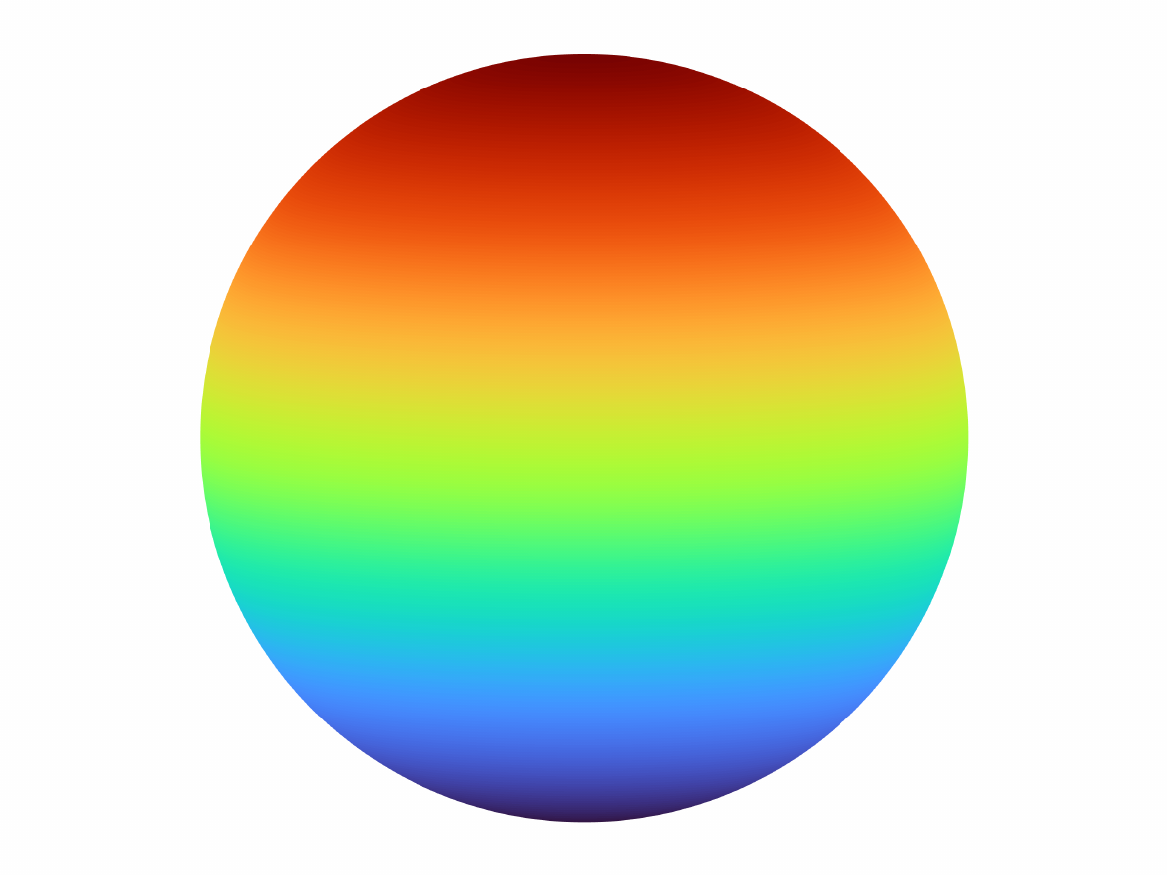}
  \end{subfigure}
  \hspace{0.5cm}
  \begin{subfigure}[b]{0.48\textwidth}
    \includegraphics[width=\textwidth]{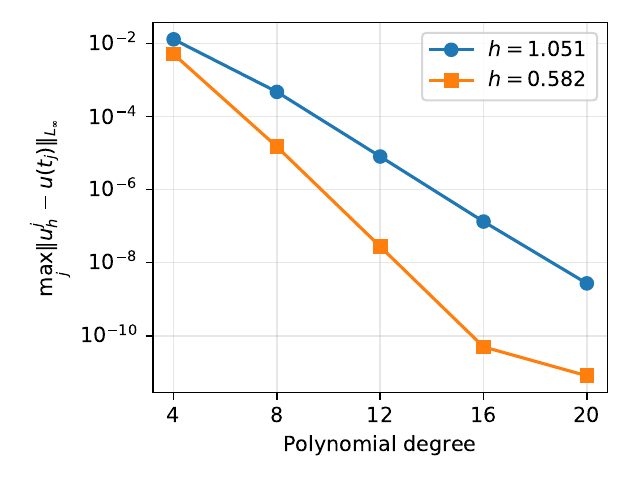}
  \end{subfigure}
  \caption{Surface diffusion on the sphere. Left: computed solution at $t=1$. Right: $E_{\max}=\max_{0\leq j\leq1000}\lVert u_h^j-u(t_j)\rVert_{L_\infty}$ versus polynomial degree for two mesh sizes.}
  \label{fig:diffsphere_error}
\end{figure}
Since the solution decays in time, Figure~\ref{fig:diffsphere_error} reports the largest nodal error over all time levels rather than the error at the final time; in all runs the maximum is attained near $t\approx0.5$, before the solution has decayed. It decreases rapidly with the polynomial degree on both meshes and saturates at the temporal error of the fixed time step.

\paragraph{Numerical study of spatial pattern
formation in Turing systems}\label{turing_model}

Having established the accuracy of the method for surface diffusion, we now turn to the Turing model of spatial pattern formation~\cite{turing1990chemical}. Turing showed that two reacting and diffusing chemicals can form spatial patterns from a nearly homogeneous initial state; an overview of common reaction kinetics is given in~\cite{maini1997spatial}. We consider the activator--inhibitor system of~\cite{barrio1999two},
\begin{equation}\label{eq:reaction_diffusion}
\begin{aligned}
    \frac{\partial u_1}{\partial t} &= \delta_{u_1} \Delta_\Gamma u_1 + \alpha u_1 \left(1 - r_1 u_2^2\right) + u_2 \left(1 - r_2 u_1\right), \\
    \frac{\partial u_2}{\partial t} &= \delta_{u_2} \Delta_\Gamma u_2 + \beta u_2 \left(1 + \frac{\alpha r_1}{\beta} u_1 u_2\right) + u_1 \left(\gamma + r_2 u_2\right),
\end{aligned}
\end{equation}
with initial data $u_1^0$ and $u_2^0$. This is~\eqref{eq:main_diff} with $\mathbf u=(u_1,u_2)$, $\mathbf D=\operatorname{diag}(\delta_{u_1},\delta_{u_2})$, and $\mathbf F$ given by the reaction terms. Here $u_1$ is the activator and $u_2$ the inhibitor, and for $\alpha=-\gamma$ the state $(u_1,u_2)=(0,0)$ is the unique equilibrium. The reaction terms create concentration peaks and diffusion smooths them; Turing patterns form when
\begin{enumerate}[label=(\roman*)]
\item without diffusion, the system tends to a uniform, linearly stable steady state, and
\item with diffusion, this state becomes unstable to small inhomogeneous perturbations.
\end{enumerate}
The patterns depend on the geometry~\cite{bunow1980pattern} and on the quadratic and cubic reaction terms, so we solve~\eqref{eq:reaction_diffusion} on several curved surfaces; a more detailed discussion of the patterns is given in~\cite{zavalani2025fast}. The cubic parameter $r_1$ promotes stripes and the quadratic parameter $r_2$ favors spots~\cite{barrio1999two}.

\begin{figure}[h!]
    \centering
    \setlength{\tabcolsep}{2pt}
    \renewcommand{\arraystretch}{0.5}
    \begin{tabular}{ccc}
     \includegraphics[width=0.32\textwidth]{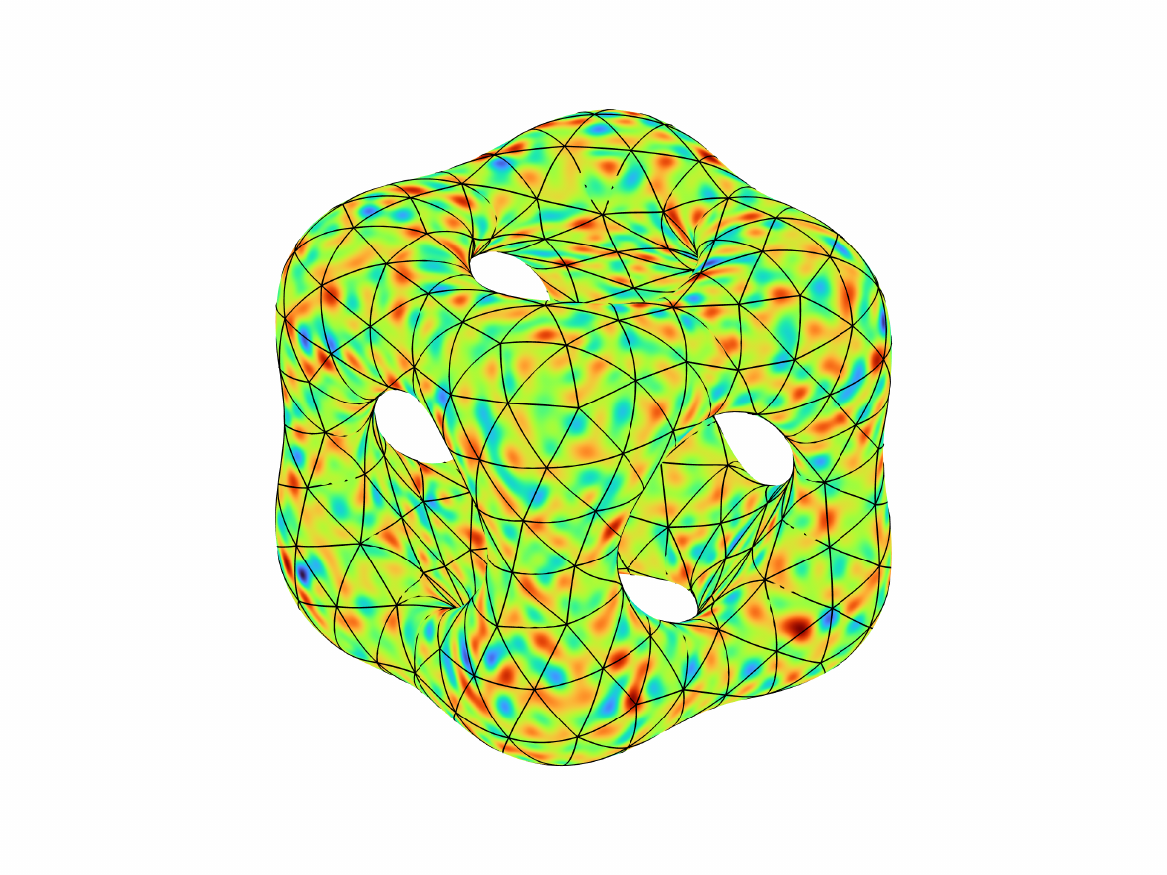} &
        \includegraphics[width=0.32\textwidth]{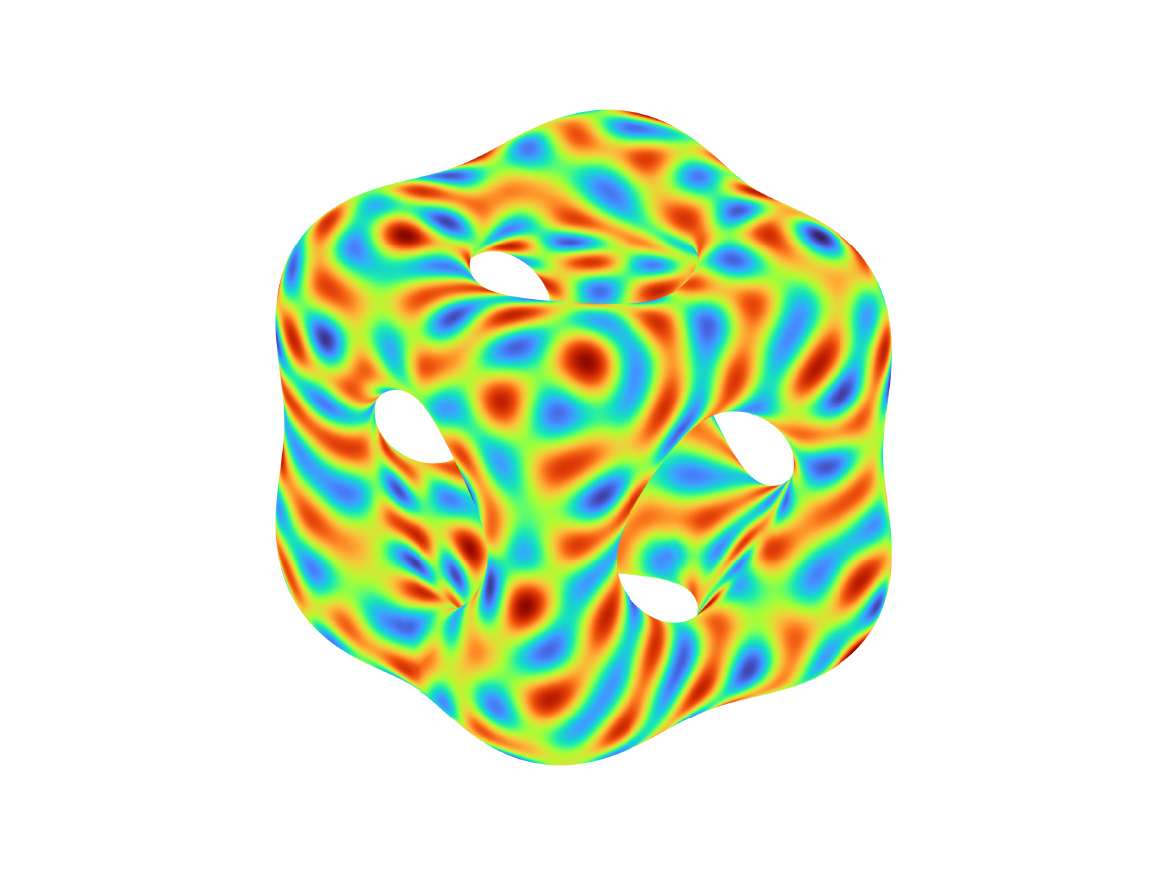} &
        \includegraphics[width=0.32\textwidth]{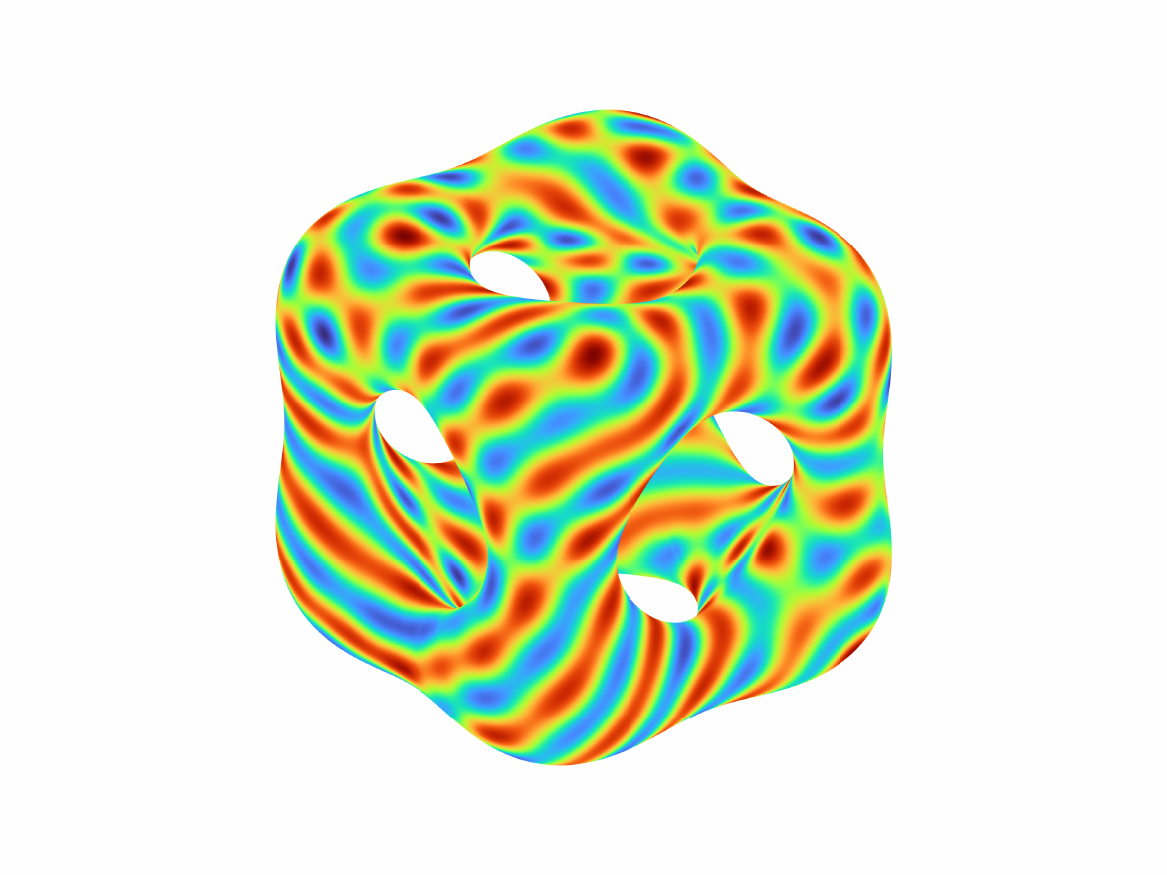} \\
        \includegraphics[width=0.32\textwidth]{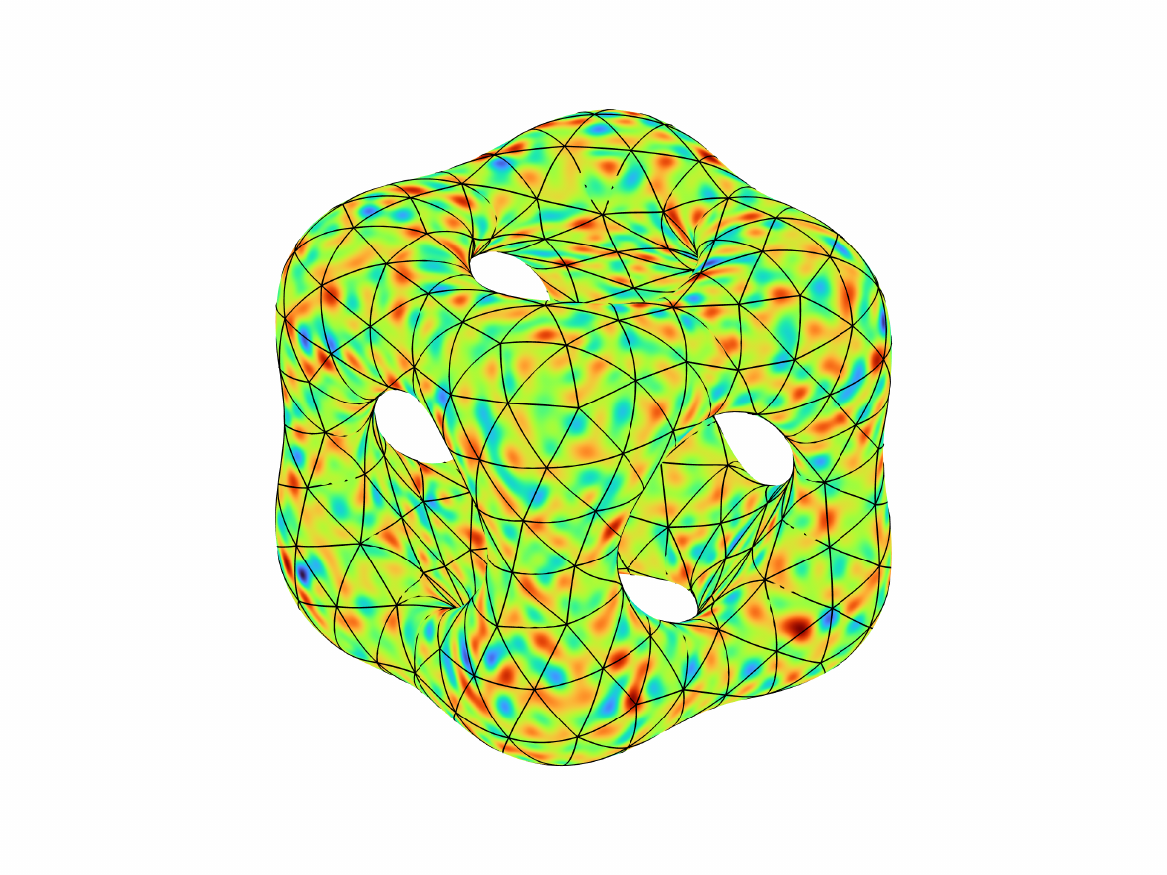} &
        \includegraphics[width=0.32\textwidth]{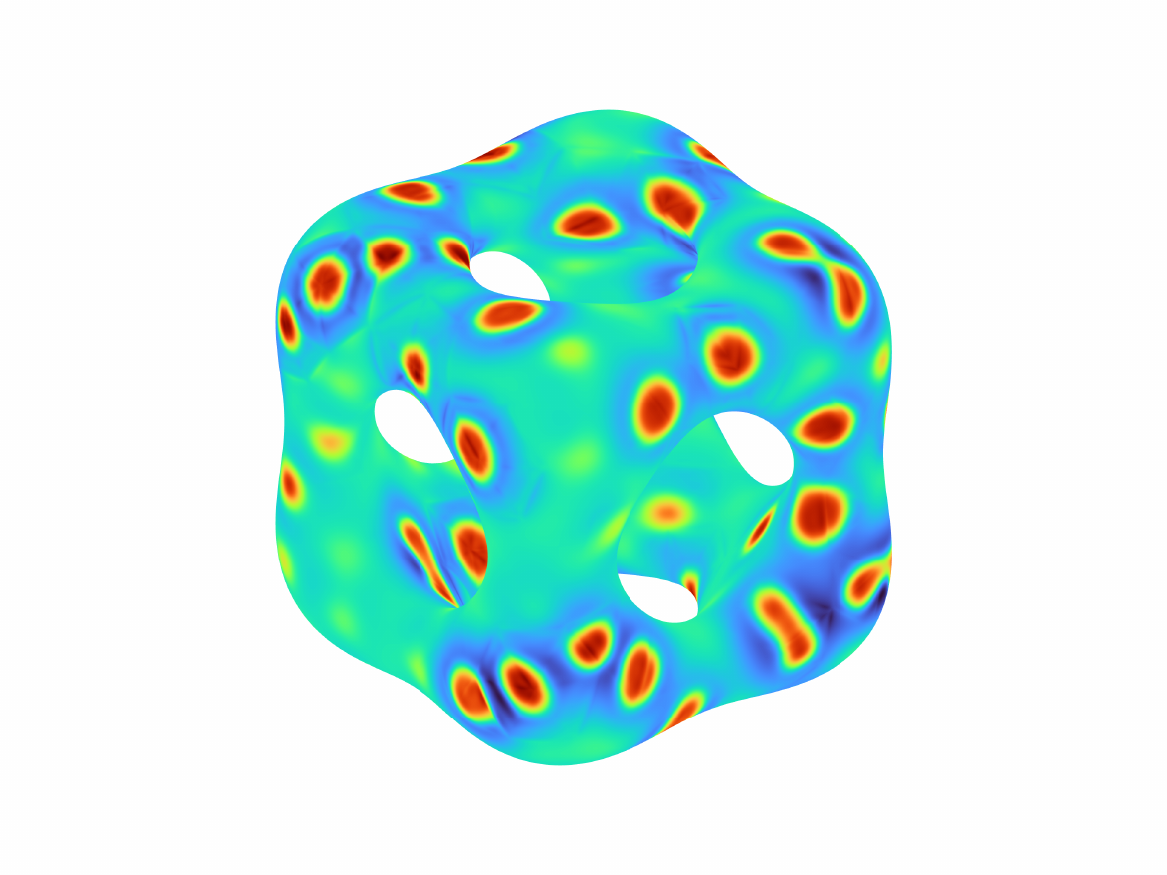} &
        \includegraphics[width=0.32\textwidth]{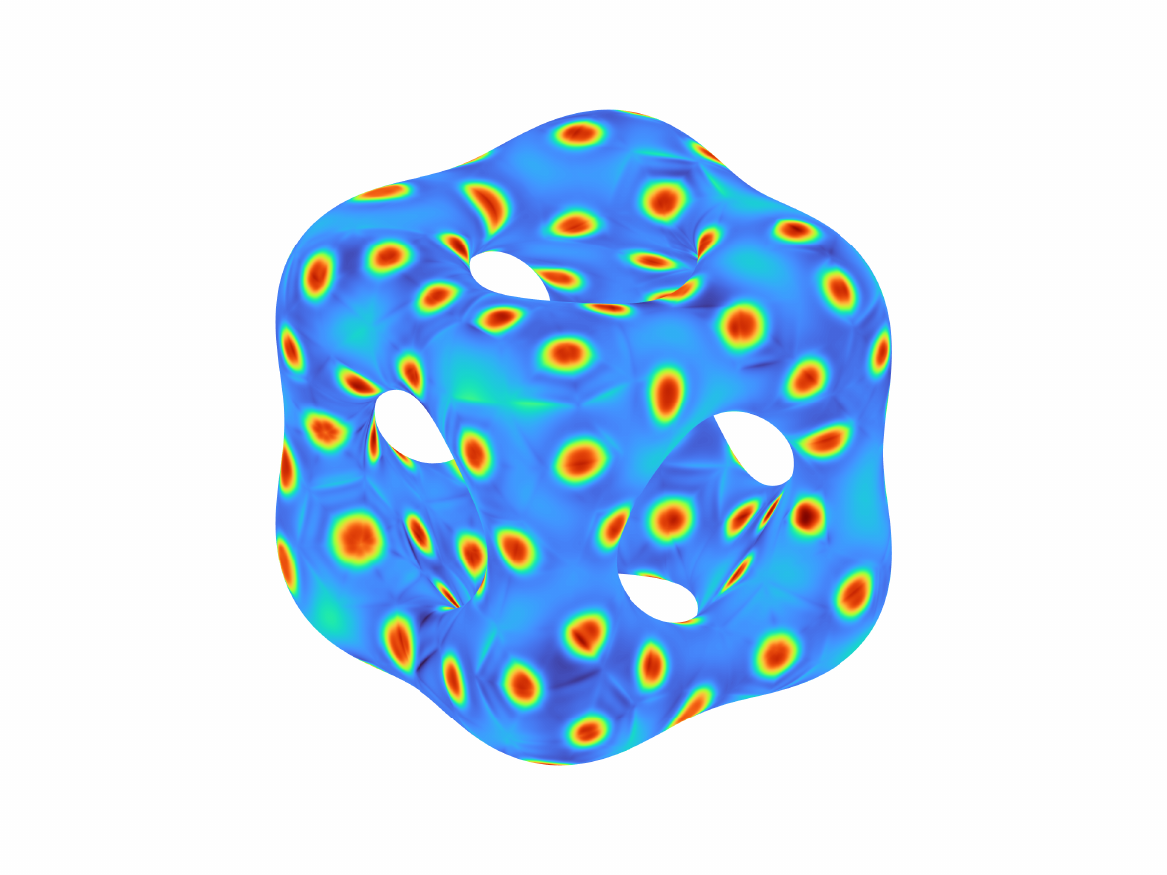}
    \end{tabular}
    \vspace{-0.8em}
    \caption{Turing model solution \( u_1 \) on a Swiss cheese surface. Top: \( r_1 = 3.5 \), \( r_2 = 0 \) at \( t = 0, 200, 600 \). Bottom: \( r_1 = 0.02 \), \( r_2 = 0.2 \) at \( t = 0, 20, 200 \). Simulated with IMEX-BDF4.}
    \label{fig:Turing_system_swiss}
\end{figure}

\begin{figure}[htbp]
  \centering
  \includegraphics[width=0.6\textwidth]{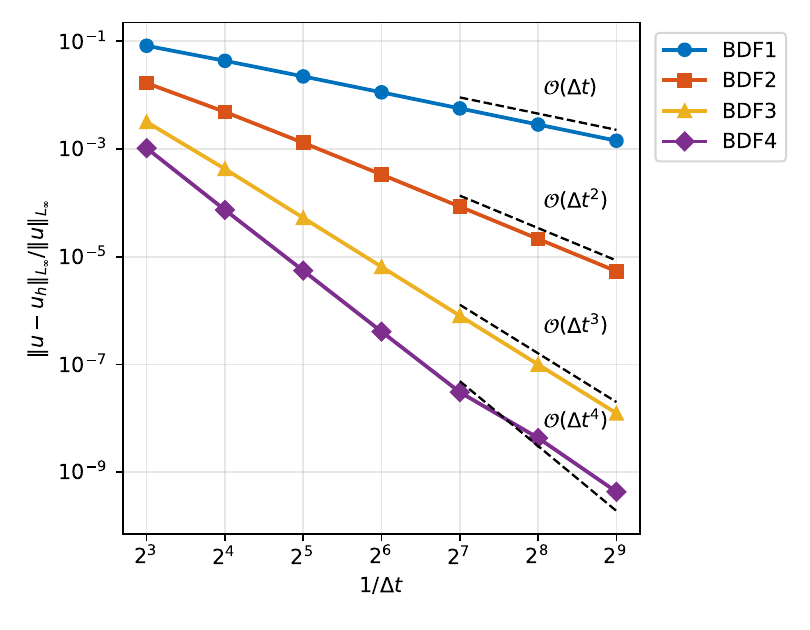}
  \caption{Temporal self-convergence of IMEX-BDF1--4 for the Turing system~\eqref{eq:reaction_diffusion} on the Swiss cheese surface ($r_1=0.02$, $r_2=0.2$). The reference solution is computed with IMEX-BDF4 and $\Delta t=2^{-10}$ to a final time $t=1$.}
  \label{fig:turing_selfconv}
\end{figure}

In our simulations, spots are more robust than stripes and reach a steady state much faster, whereas stripes emerge only for very small $r_2$ and their orientation depends on the initial data. \Cref{fig:Turing_system_swiss} shows both regimes on a Swiss cheese surface\footnote{The Swiss cheese surface is defined implicitly by $f(\mathbf{x}) = (x^2 + y^2 - 4)^2 + (z^2 - 1)^2 + (y^2 + z^2 - 4)^2 + (x^2 - 1)^2 + (z^2 + x^2 - 4)^2 + (y^2 - 1)^2 - 15$.} for the dimensionless parameters of~\cite{bunow1980pattern}, $\alpha=0.899$, $\beta=-0.91$, $\gamma=-\alpha$, $\delta_{u_2}=5\cdot10^{-3}$, and $\delta_{u_1}=0.516\,\delta_{u_2}$, with random initial data.
For the spot regime, we also check that the time integration is accurate and not only stable by comparing IMEX-BDF1--4 with a reference solution (\Cref{fig:turing_selfconv}). IMEX-BDF$k$ converges at order $k$, so the explicit treatment of the reaction terms does not reduce the order, and IMEX-BDF4 is fourth-order accurate with one factorization per species.

\begin{figure}[htbp]
	\centering
        \setlength{\tabcolsep}{2pt}
    \renewcommand{\arraystretch}{0.5}
	\begin{tabular}{ccc}
        	~~\includegraphics[width=0.31\textwidth]{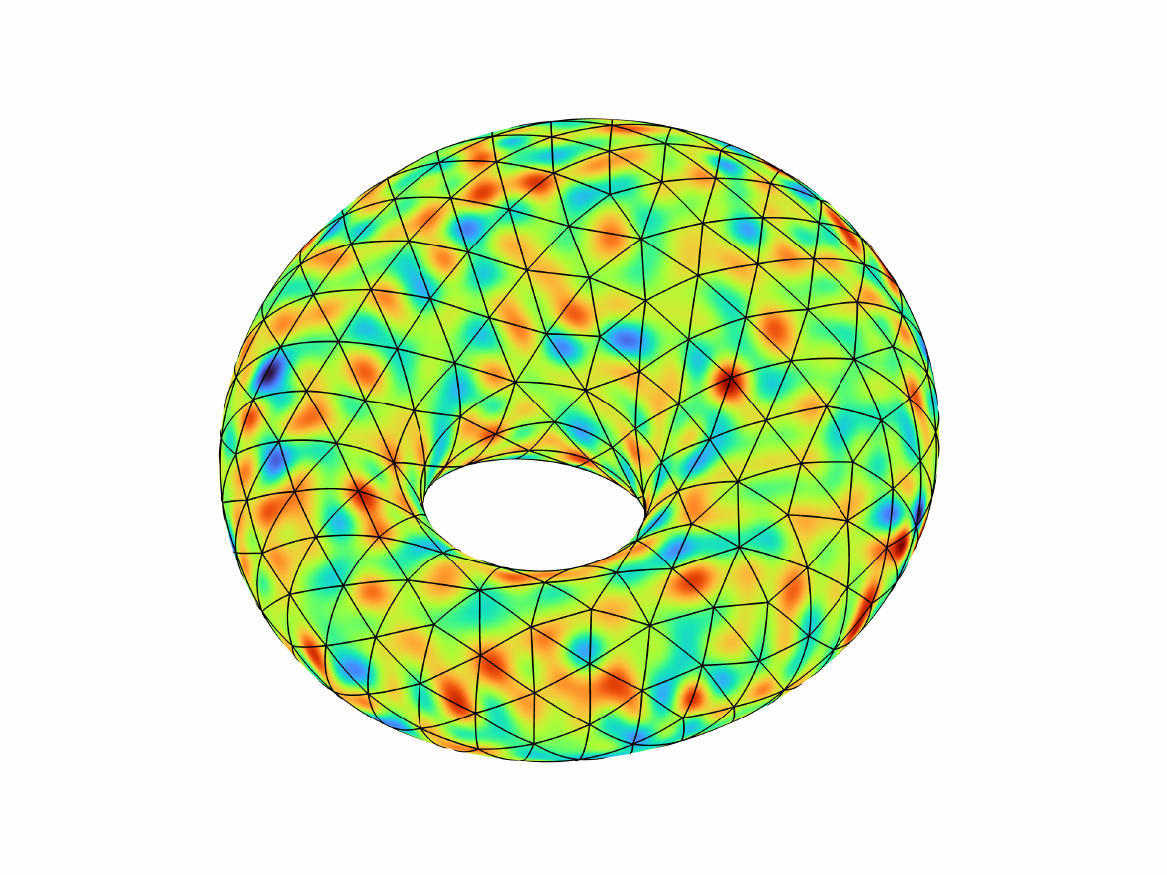} &
	~\includegraphics[width=0.31\textwidth]{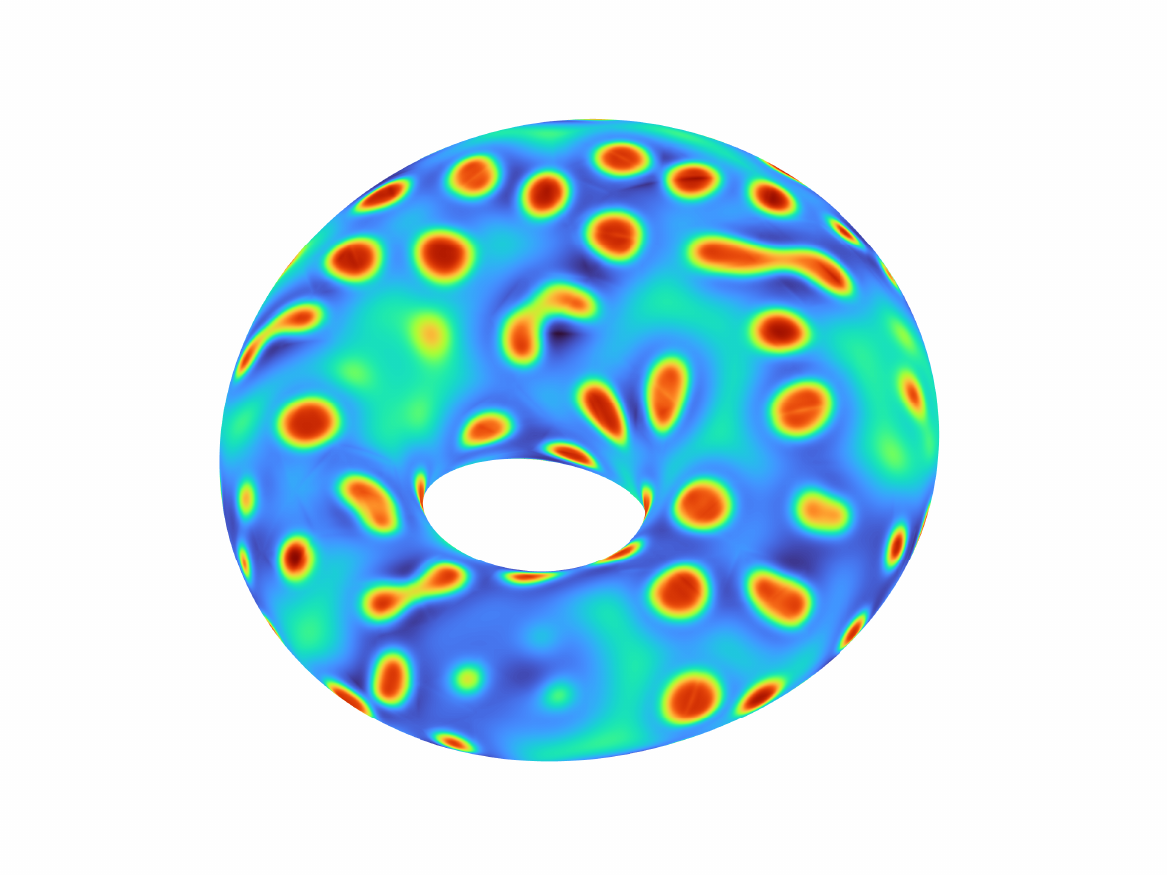} &
	~~\includegraphics[width=0.31\textwidth]{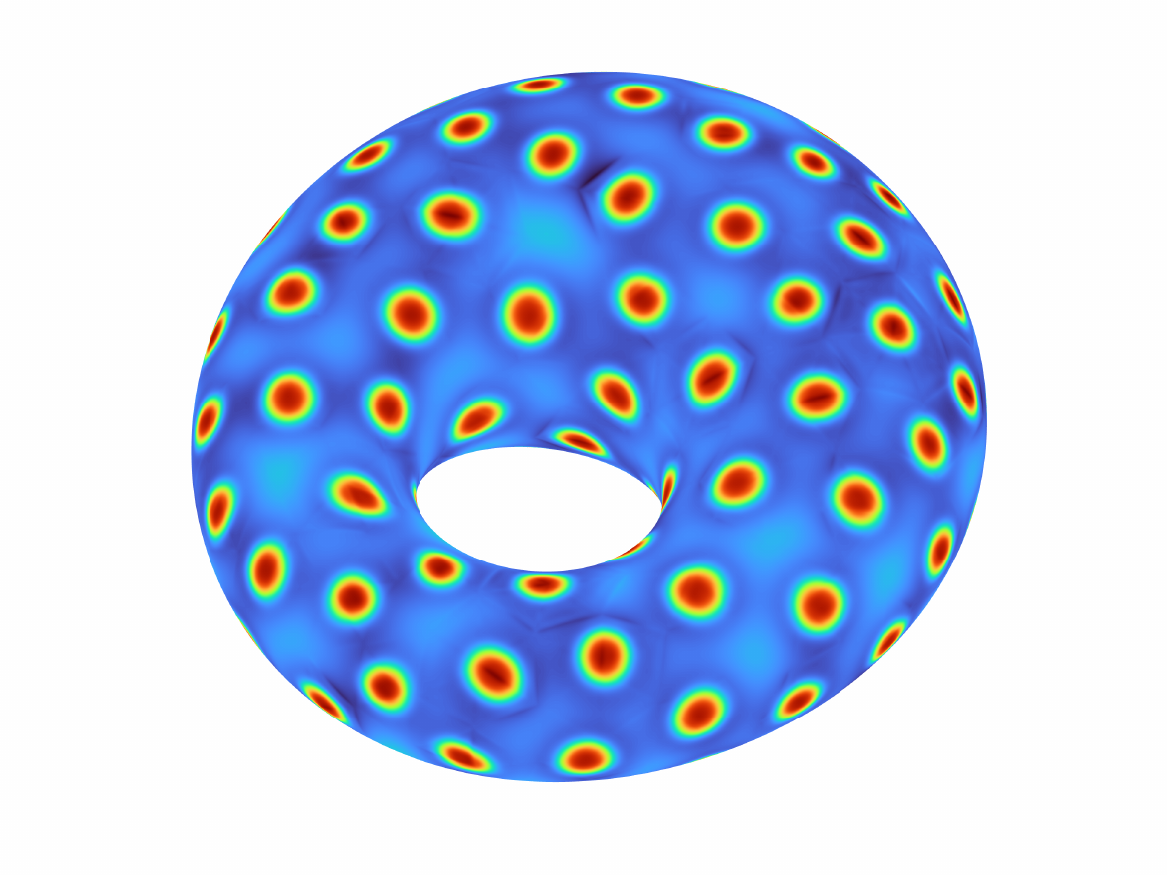}
    \\
	
    \includegraphics[width=0.31\textwidth]{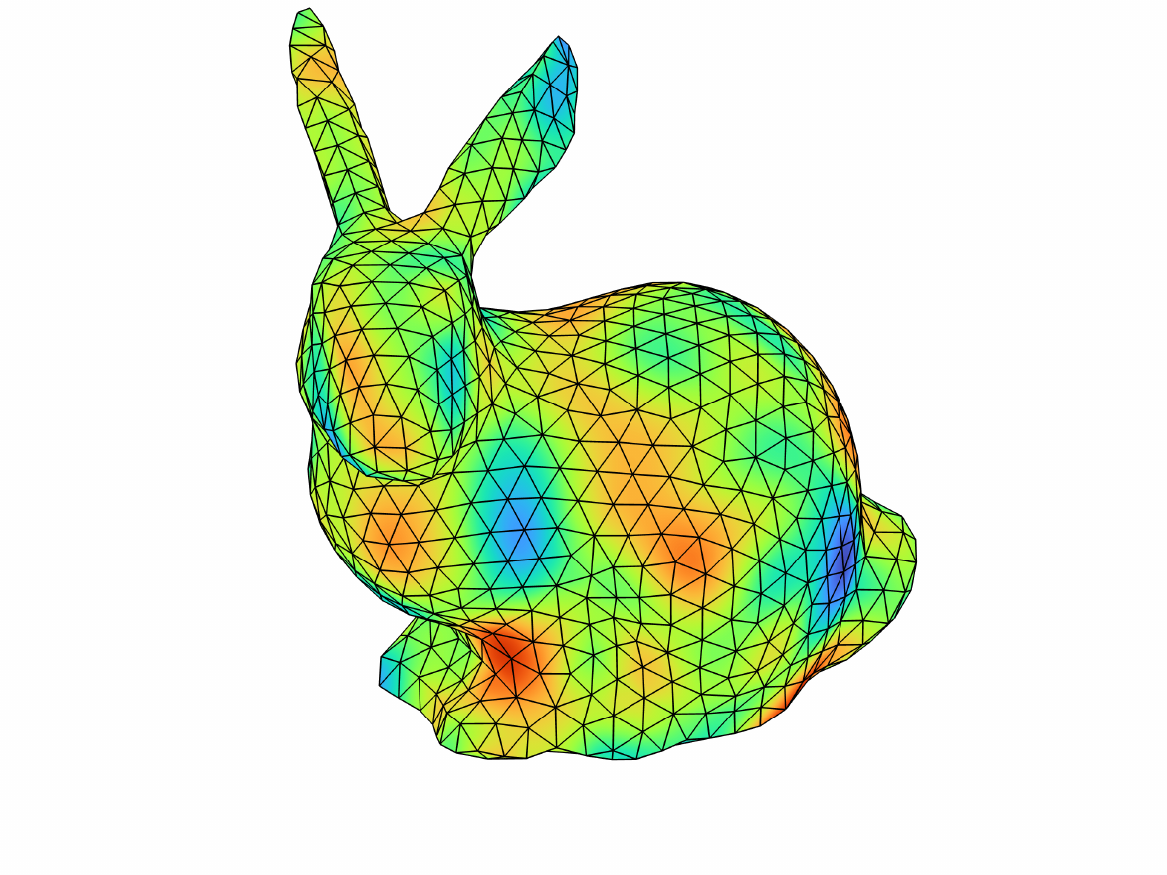} &
	\includegraphics[width=0.31\textwidth]{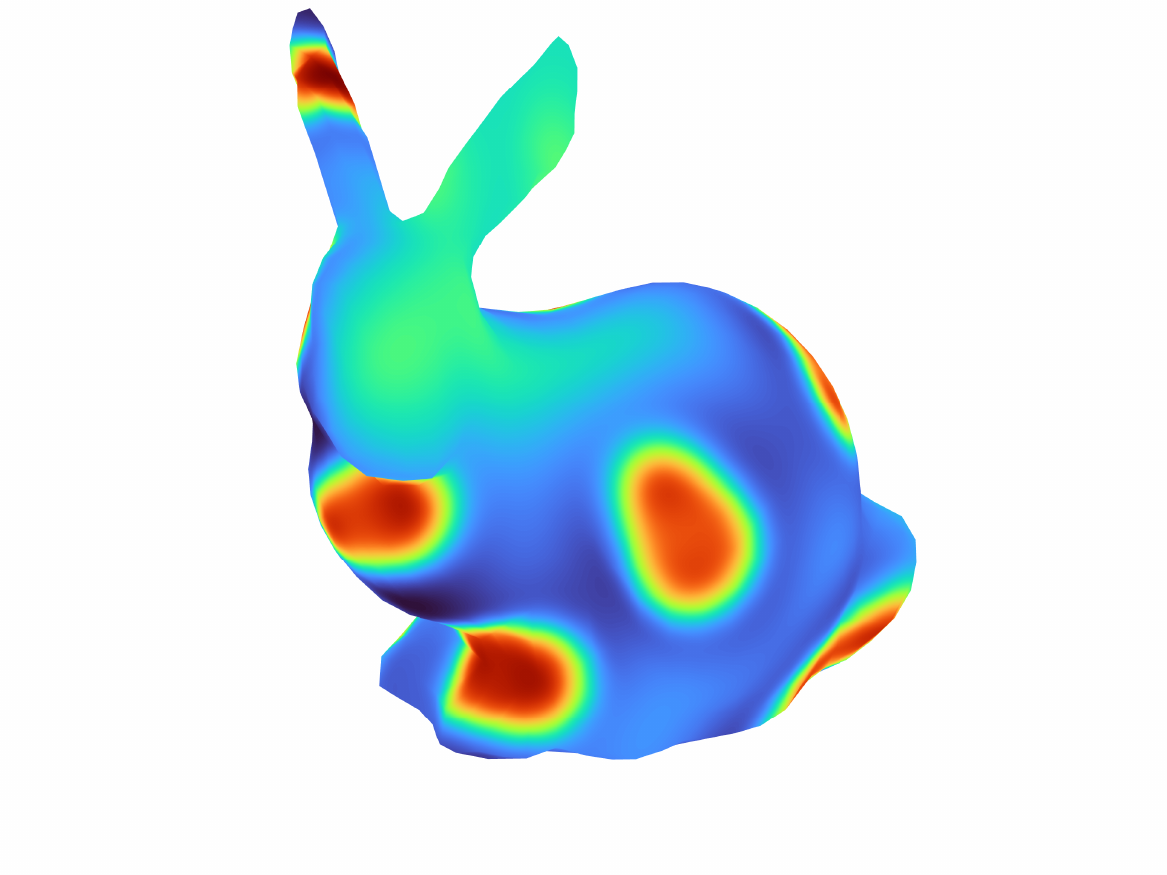} &
	\includegraphics[width=0.31\textwidth]{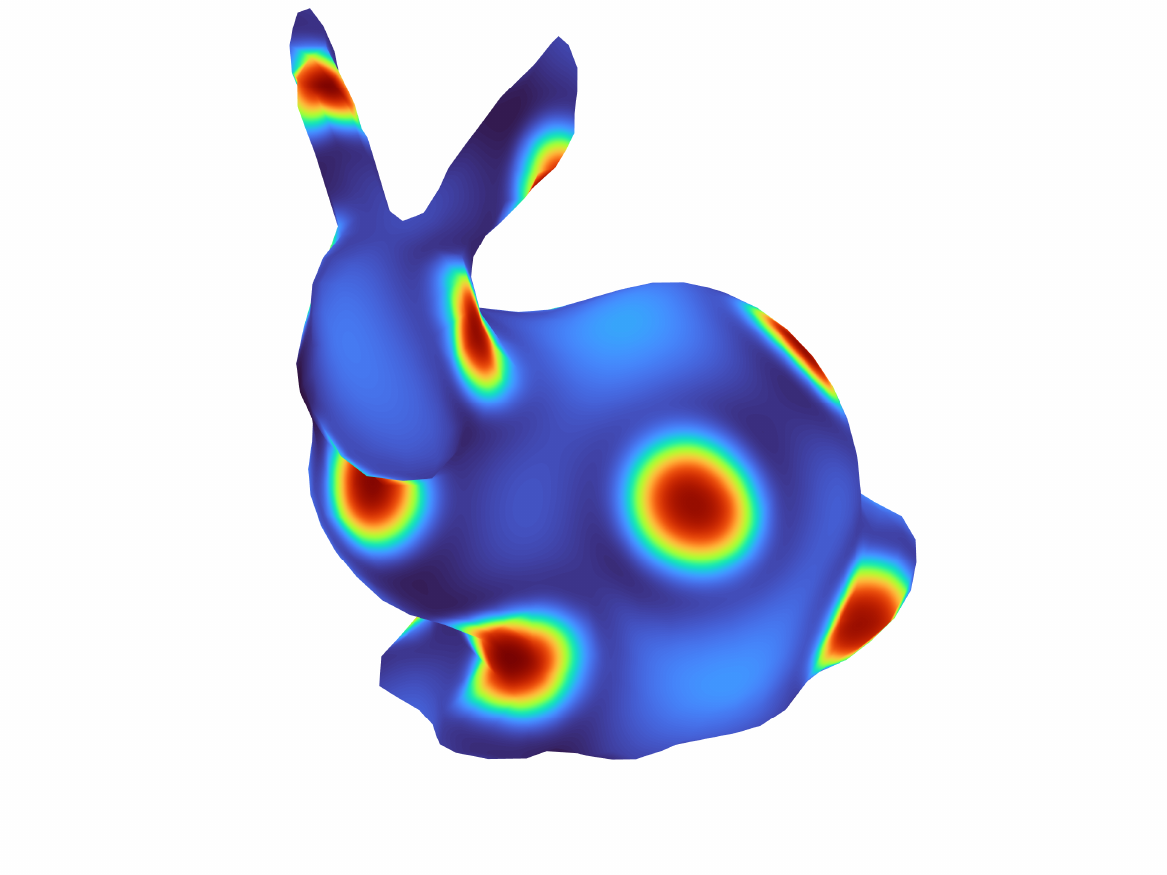} 

	\end{tabular}
\caption{Turing model solution \( u_1 \) on the asymmetric torus and Stanford Bunny at times \( t = 0 \), \( 20 \), and \( 200 \) using IMEX-BDF4 with \( r_1 = 0.02 \), \( r_2 = 0.2 \).}

\label{fig:Turing system_combined}
\end{figure}

Next, we study the influence of the geometry by solving the Turing model on an asymmetric torus\footnote{The surface is given implicitly by $f(\mathbf{x}) = \bigl(x^2 + y^2 + z^2 - d^2 + b^2\bigr)^2 - 4\bigl(a x + c^2 d\bigr)^2 - 4 b^2 y^2$, with $a = 2$, $b = 1.9$, $d = 1$, and $c^2 = a^2 - b^2$.} and on a triangulated Stanford Bunny with the same parameters, $r_1=0.02$ and $r_2=0.2$. We use elements of degree~$10$ and IMEX-BDF4 with $\Delta t=0.1$ for $2000$ steps, up to $t=200$. In \Cref{fig:Turing system_combined}, the rows show the two surfaces and the columns $t=0$, $20$, and $200$. Although the parameters are identical, the spots emerge and arrange themselves differently on the two surfaces, so the pattern also depends on the shape of the surface. This dependence is often missed in one- and two-dimensional studies, where the geometry is simplified or ignored.

We conclude this paragraph with interacting Turing systems. To our knowledge, such systems have so far been simulated only on planar domains~\cite{barrio1999two}. Following~\cite{barrio1999two}, a second Turing system in the chemicals $(u_1,u_2)$ modulates the kinetic terms of the $(v_1,v_2)$ system, giving the model:
\begin{equation}\label{eq:interaction}
\begin{aligned}
\frac{\partial v_1}{\partial t} &= \delta_{v_1} \Delta_{\Gamma} v_1 + \alpha' v_1 (1 - r_1 v_2^2) + v_2(1 - r_2 v_1) + q_1 u_1 + q_2 u_1 v_2 + q_3 u_1 v_2^2 \\
\frac{\partial v_2}{\partial t} &= \delta_{v_2} \Delta_{\Gamma} v_2 + \beta' v_2 \left( 1 + \frac{\alpha' r_1}{\beta'} v_1 v_2 \right) + v_1(\gamma' + r_2 v_2) - q_2 u_2 v_1 - q_3 u_2^2 v_1 \\
\frac{\partial u_1}{\partial t} &= \delta_{u_1} \Delta_{\Gamma} u_1 + \alpha u_1(1 - r_1 u_2^2) + u_2(1 - r_2 u_1) \\
\frac{\partial u_2}{\partial t} &= \delta_{u_2} \Delta_{\Gamma} u_2 + \beta u_2 \left( 1 + \frac{\alpha r_1}{\beta} u_1 u_2 \right) + u_1(\gamma + r_2 u_2).
\end{aligned}
\end{equation}
This setup allows us to investigate how pattern properties change due to the interaction between the two coupled systems. For the simulations, we use the following parameter values:

\begin{itemize}
  \item For the \( (v_1, v_2) \) system:
  \[
  \alpha' = 0.398, \quad \beta' = -0.41, \quad \gamma' = -\alpha', \quad 
  \delta_{v_2} = 5 \times 10^{-3}, \quad \delta_{v_1} = 0.122\, \delta_{v_2}.
  \]
  
  \item For the \( (u_1, u_2) \) system:
  \[
  \alpha = 0.899, \quad \beta = -0.91, \quad \gamma = -\alpha, \quad 
  \delta_{u_2} =  \delta_{v_2}, \quad \delta_{u_1} = 0.516\, \delta_{u_2}.
  \]
\end{itemize}

\begin{figure*}[!t]
  \centering

  \begin{subfigure}{0.24\textwidth}
    \centering
    \includegraphics[width=\linewidth]{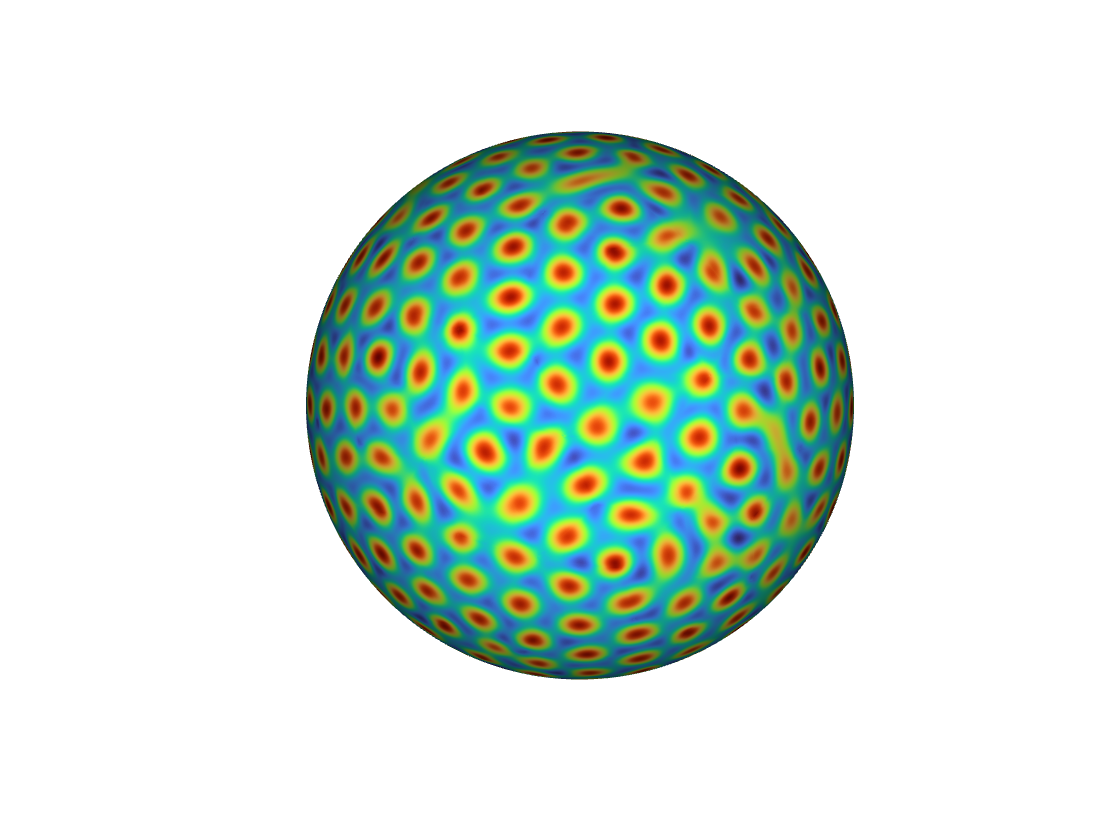}
    \caption{ $u_1$}
  \end{subfigure}\hfill
  \begin{subfigure}{0.24\textwidth}
    \centering
    \includegraphics[width=\linewidth]{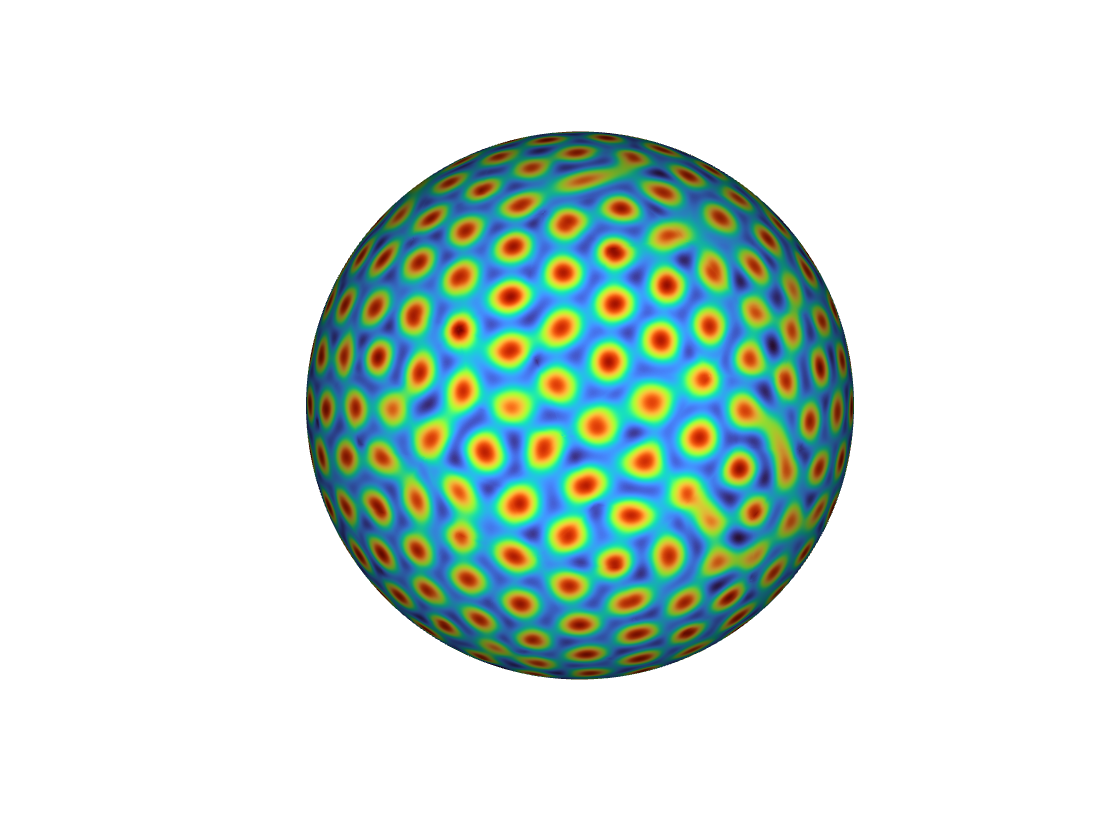}
    \caption{$v_1$, with $q_1 \neq 0$}
  \end{subfigure}\hfill
  \begin{subfigure}{0.24\textwidth}
    \centering
    \includegraphics[width=\linewidth]{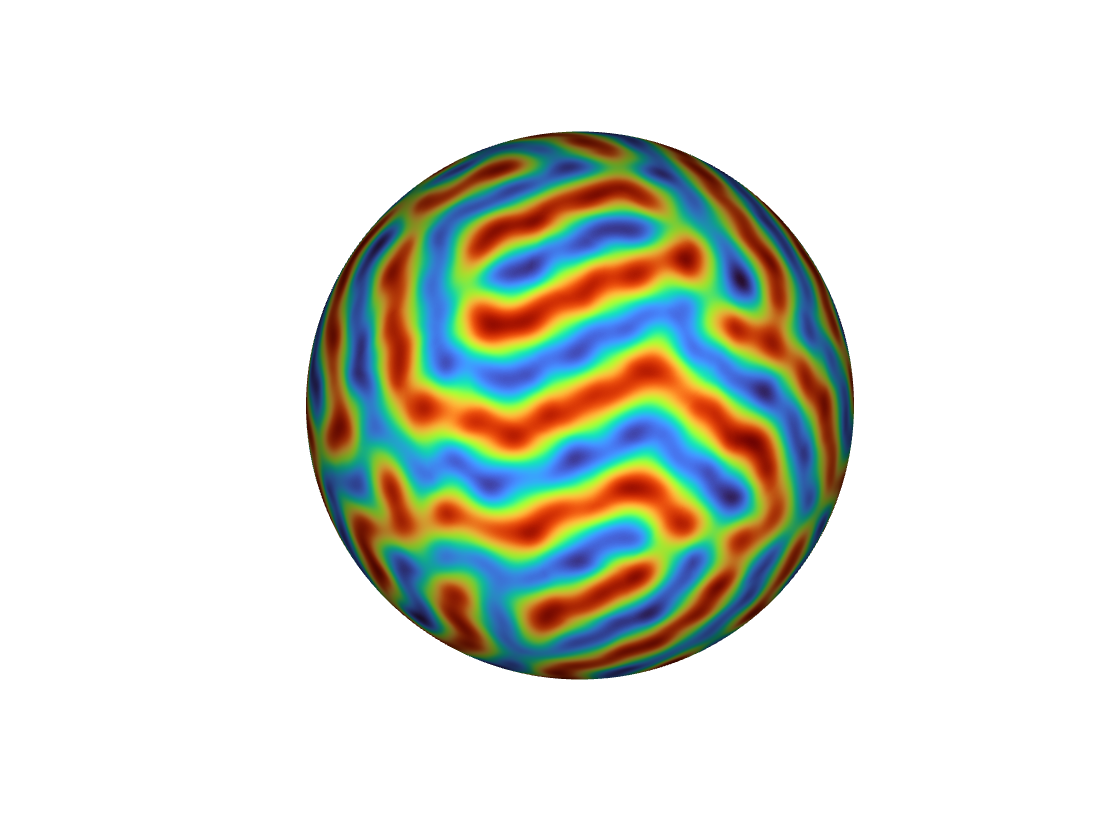}
    \caption{$v_1$, with $q_2 \neq 0$}
  \end{subfigure}\hfill
  \begin{subfigure}{0.24\textwidth}
    \centering
    \includegraphics[width=\linewidth]{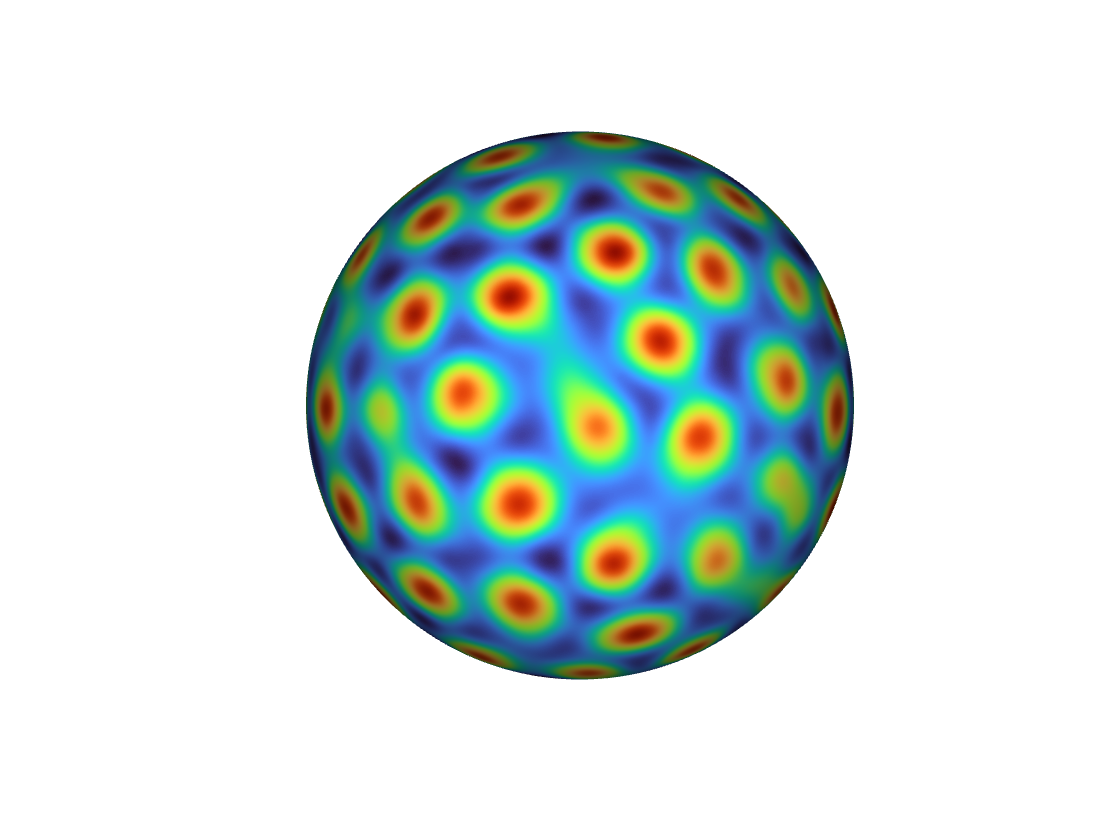}
    \caption{$v_1$, with $q_3 \neq 0$}
  \end{subfigure}

    \caption{
    Patterns obtained from the coupled system~\eqref{eq:interaction}.
    Each panel isolates a single coupling mechanism:
    linear ($q_1$), quadratic ($q_2$), or cubic ($q_3$),
    with coupling strength $0.55$.
  }
 \label{fig:coupled_patterns}
\end{figure*}

\begin{figure*}[!t]
  \centering

  \begin{subfigure}{0.24\textwidth}
    \centering
    \includegraphics[width=\linewidth]{images/up_1_q.png}
    \caption{$u_1$}
  \end{subfigure}\hfill
  \begin{subfigure}{0.24\textwidth}
    \centering
    \includegraphics[width=\linewidth]{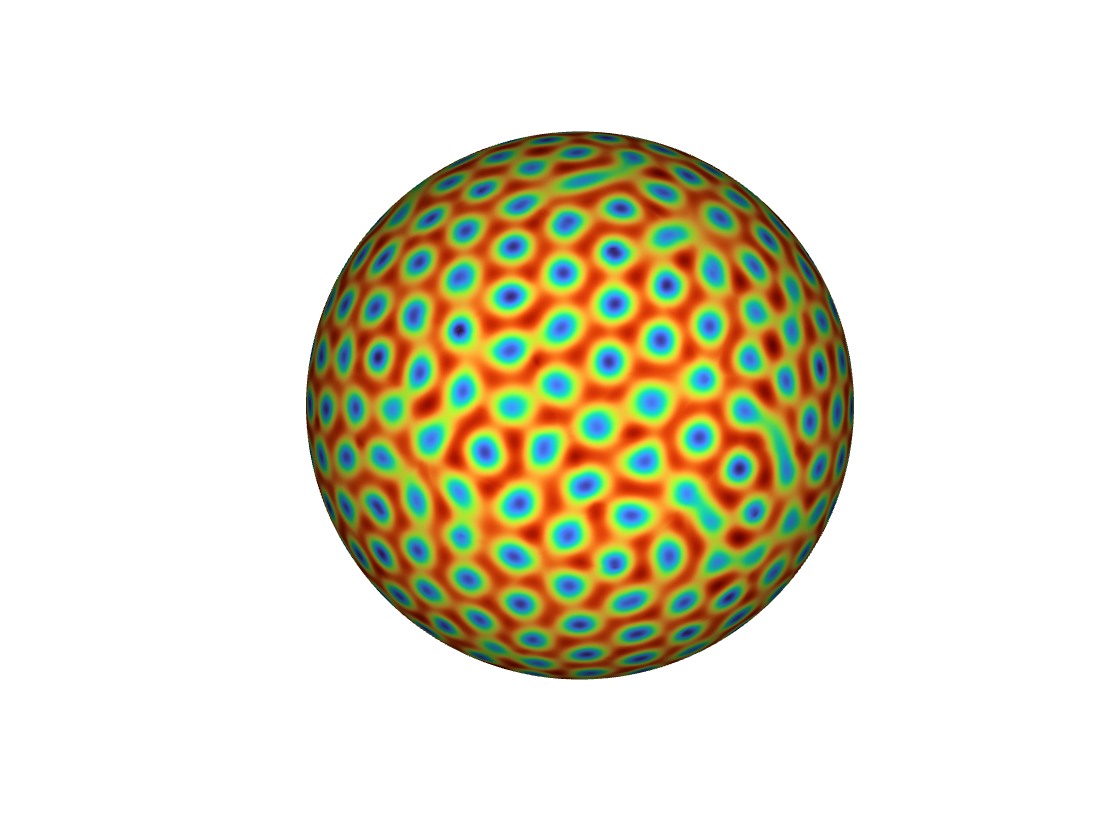}
    \caption{$v_1$, with $q_1 \neq 0$}
  \end{subfigure}\hfill
  \begin{subfigure}{0.24\textwidth}
    \centering
    \includegraphics[width=\linewidth]{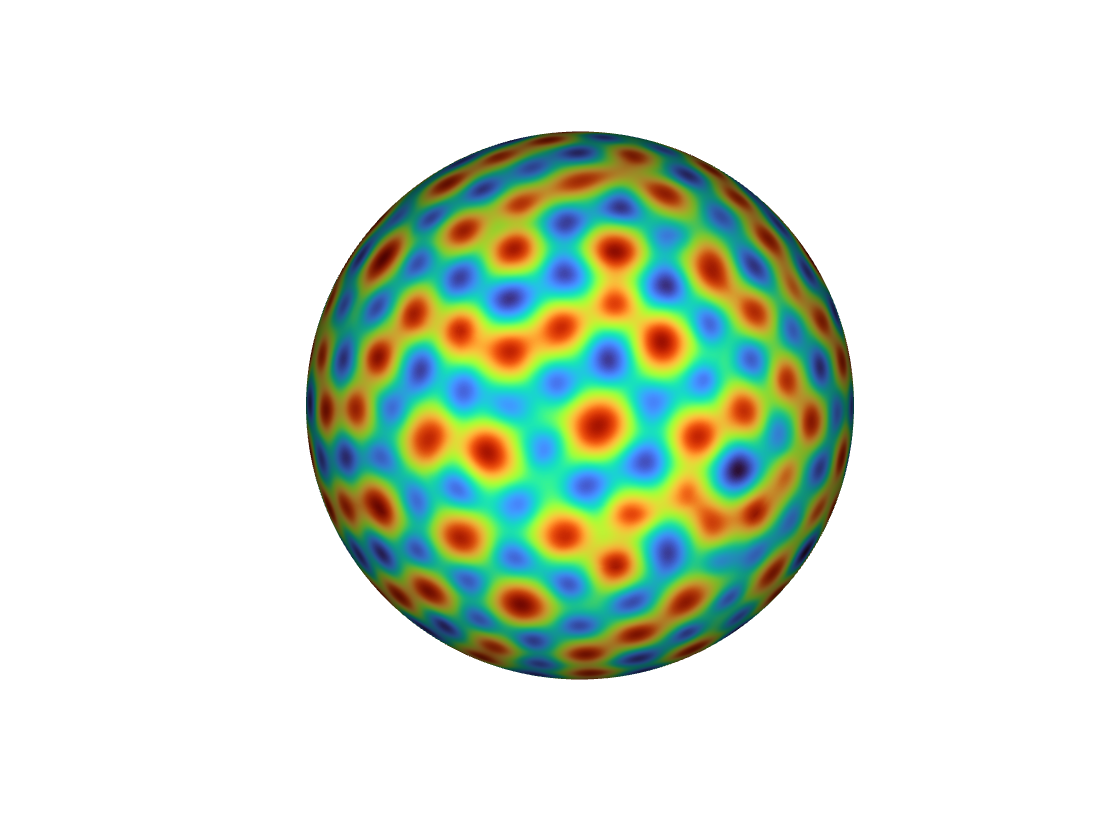}
    \caption{$v_1$, with $q_2 \neq 0$}
  \end{subfigure}\hfill
  \begin{subfigure}{0.24\textwidth}
    \centering
    \includegraphics[width=\linewidth]{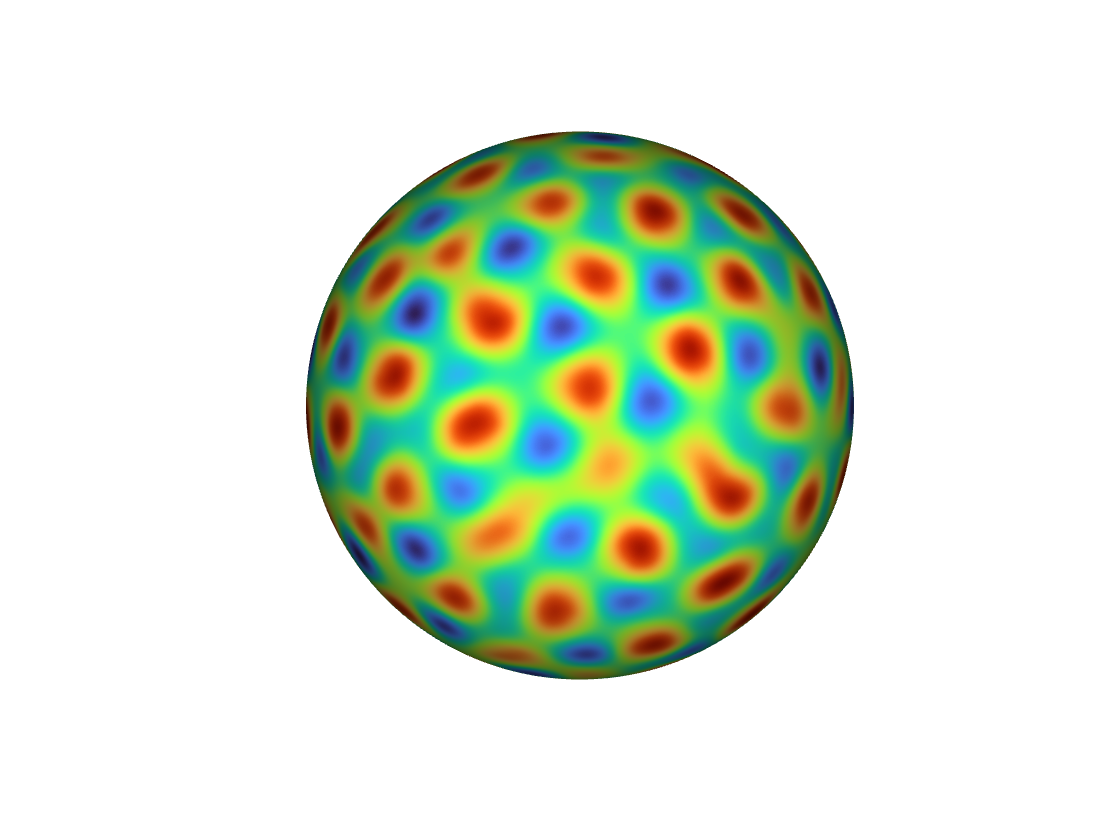}
    \caption{$v_1$, with $q_3 \neq 0$}
  \end{subfigure}

    \caption{
    Patterns obtained from the coupled system~\eqref{eq:interaction},
    using the same parameters as in Fig.~\ref{fig:coupled_patterns},
    but with negative coupling coefficients set to $-0.55$.
  }
\label{fig:-coupled_patterns}
\end{figure*}

\Cref{fig:coupled_patterns} shows the effect of each coupling. Linear coupling makes the pattern of $v_1$ identical to that of $u_1$, so the coupling overrides the dynamics of $v_1$. Cubic coupling leaves the spot pattern nearly unchanged. Quadratic coupling has the strongest effect: it distorts the spots and overlays them on a background of stripes, giving a mixed pattern of spots and labyrinths. With negative coupling coefficients (\Cref{fig:-coupled_patterns}), linear coupling produces an inverted image of $u_1$ in $v_1$. The other two couplings superpose stripes and spots: quadratic coupling now favors spots, and cubic coupling arranges them more regularly.

\section{Conclusions and future work}\label{sec:conclusion}

We have presented a high-order fast direct solver for PDEs on triangulated surfaces. On each curved triangle the operator is discretized by collocation in the Dubiner basis, and the elements are coupled on a vertex-free edge grid and merged hierarchically, so that a general triangulation can be used without first converting it into quadrilaterals. We observe convergence at the rate $\mathcal O(h^{n-2})$ under mesh refinement. The factorization costs $\mathcal O(N_{\mathrm{dof}}^{3/2})$ at fixed degree and is then reused for every right-hand side at a cost of $\mathcal O(N_{\mathrm{dof}}\log N_{\mathrm{dof}})$, which makes the solver well suited to implicit time stepping on complex geometries.

We plan to extend the solver to vector-valued problems, such as the surface Stokes equations, and to evolving surfaces, where the element maps, and with them the local solution and Dirichlet-to-Neumann operators, change at every time step.

\section{Appendix. Tangential differential calculus on surfaces}\label{tang_operator}

We review the basic concepts of differential geometry relevant to our setting. We start by describing the parametric representation, then  introduce the first fundamental form $\mathbf{g}$, the tangential operators (gradient $\nabla_\Gamma$, divergence $\operatorname{div}_{\Gamma} $, and Laplace–Beltrami $\Delta_\Gamma$). Let $\Gamma$ be a smooth surface embedded
in $\R^{d+1}$ locally parametrized by the mapping 
\begin{equation*}
    \mathbf{X}_{i}: \Delta_{d}\longrightarrow \R^{d+1}\;, i=1,2,\ldots, N.
\end{equation*}

 Let \(\partial_j \mathbf{X}_{i} (s,t)\) be the column vector of \(j\)-th partial derivatives of \(\mathbf{X}_{i}(\mathbf{y})\) for \(1 \leq j \leq d\) at \(\mathbf{y} \in \Delta_d\).
 For immersed surfaces $\Gamma \subset \R^{d+1}$ we will write
$D\mathbf{X}_{i}(\mathbf{y}) : \R^d \to \R^{d+1}$ for the
Jacobian of the parametrization $\varrho_i$ at $\mathbf{y}$.  The rank of \(D\mathbf{X}_{i}(\mathbf{y}) = (\partial_j \mathbf{X}_{i} (\mathbf{y}))_{j=1}^{d} \in \mathbb{R}^{(d+1) \times d}\) is \(d\) (full rank). This implies that \(\{ \partial_j \mathbf{X}_{i} (\mathbf{y}) \}_{j=1}^{d}\) are linearly independent and span the tangent hyperplane to \(\Gamma\) at \(\mathbf{x}\in  \mathcal{U} \cap \Gamma \), where $\mathcal{U}$ is an open set.

The \textit{first fundamental form} is the symmetric and positive definite matrix 
\(\mathbf{g} \in \mathbb{R}^{d \times d}\) defined by
$\mathbf{g}(\mathbf{y}) := D\mathbf{X}_{i}(\mathbf{y}) ^T D\mathbf{X}_{i}(\mathbf{y}).$
If \(\mathbf{g} = (g_{ij})_{i,j=1}^{d}\), then the components \(g_{ij}\) read
\[
g_{ij} = \partial_i \mathbf{X}_{i}^T \partial_j \mathbf{X}_{i} = \partial_i \mathbf{X}_{i} \cdot \partial_j \mathbf{X}_{i},
\]
which depends on the choice of parametrization.

Given $v\in C^1(\Delta_d,\mathbb{R})$, let $\tilde v:\Gamma\to\mathbb{R}$ be the
corresponding function defined by $\tilde v := v\circ \mathbf{X}^{-1}$. 
The \emph{tangential (or surface) gradient} of $\tilde v$ is defined as the vector
tangent to $\Gamma$ that satisfies the chain rule
\[
\nabla{v}(\mathbf{y}) = D\mathbf{X}_{i}(\mathbf{y})^T \nabla_{\Gamma} \tilde{v} (\mathbf{x}) \quad \forall \mathbf{y} \in \Delta_{d}.
\]
Since \( \nabla_{\Gamma} \tilde{v} \) is spanned by \( \{ \partial_j \mathbf{X}_{i} \}_{j=1}^d \), we get \( \nabla_{\Gamma} \tilde{v} = D\mathbf{X}_{i} w \) for some \( w \in \mathbb{R}^d \) whence \( w = \mathbf{g}^{-1} \nabla v \) and

\begin{equation*}
\nabla_{\Gamma} \tilde{v} (\mathbf{x}) = D\mathbf{X}_{i}(\mathbf{y}) \mathbf{g}(\mathbf{y})^{-1} \nabla v (\mathbf{y}) \quad \forall \mathbf{y} \in \Delta_{d}. 
\end{equation*}

If \( \tilde{v} = (\tilde{v}_i)_{i=1}^{d+1} : \Gamma \rightarrow \mathbb{R}^{d+1} \) is a vector field of class \( C^1 \), we define its tangential differential \( D_{\Gamma} \tilde{v} \in \mathbb{R}^{(d+1) \times (d+1)} \) as a matrix whose \( i \)-th row is \( (\nabla_{\Gamma} \tilde{v}_i)^T \). 
\begin{equation*}\label{eq:grad}
D_{\Gamma} \tilde{v} = 
\begin{pmatrix}
(\nabla_{\Gamma} \tilde{v}_1)^T \\
(\nabla_{\Gamma} \tilde{v}_2)^T \\
\vdots \\
(\nabla_{\Gamma} \tilde{v}_{d+1})^T 
\end{pmatrix}=\begin{pmatrix}
(\partial_{1}^{\Gamma} v_1, \ldots, \partial_{d+1}^{\Gamma} v_1)^T \\
(\partial_{1}^{\Gamma} v_2, \ldots, \partial_{d+1}^{\Gamma} v_2)^T \\
\vdots \\
(\partial_{1}^{\Gamma} v_3, \ldots, \partial_{d+1}^{\Gamma} v_3)^T
\end{pmatrix},
\end{equation*}
for the $d+1$ components of the tangential gradient. Here, we identify $\partial_i^{\Gamma}$ with the tangential derivative $\partial_i$ along the surface $\Gamma$.

In addition, the \textit{tangential divergence} of \( \tilde{v} \) is the trace of \( D_{\Gamma} \tilde{v} \)
\begin{equation*}
\operatorname{div}_{\Gamma} \tilde{v} (\mathbf{x}) = \operatorname{trace}(D_{\Gamma} \tilde{v} (\mathbf{x})) = \sum_{i,j=1}^{d} g^{ij} (\mathbf{y}) \partial_i \mathbf{X}_{i} (\mathbf{y}) \cdot \partial_j v (\mathbf{y}) \quad \forall \mathbf{y} \in \Delta_{d}. 
\end{equation*}

provided \( g^{-1} = (g^{ij})_{i,j=1}^d \). If both \( \Gamma \) and \( v : \Gamma \rightarrow \mathbb{R} \) are of class \( C^2 \), then the \textit{Laplace--Beltrami (or surface Laplace) operator} is now defined to be
\begin{equation*}
\Delta_{\Gamma} \tilde{v} = \frac{1}{\sqrt{\det \mathbf{g}(\mathbf{y})}} \operatorname{div} \left( \sqrt{\det \mathbf{g}(\mathbf{y})} \mathbf{g}(\mathbf{y})^{-1} \nabla v (\mathbf{y}) \right) \quad \forall \mathbf{y} \in \Delta_{d}. 
\end{equation*}

\section*{Funding}
This work was supported by the German Research Foundation (DFG) through the project ``Beyond uni-directional hierarchical wrinkle formation'' (project number 386450667).
\bibliographystyle{siamplain}
\bibliography{Ref.bib}
\end{document}